# Lecture notes on
# 2-dimensional defect TQFT

Nils Carqueville

nils.carqueville@univie.ac.at

*Fakultät für Mathematik, Universität Wien, Austria*

These notes offer an introduction to the functorial and algebraic description of 2-dimensional topological quantum field theories 'with defects', assuming only superficial familiarity with closed TQFTs in terms of commutative Frobenius algebras. The generalisation of this relation is a construction of pivotal 2-categories from defect TQFTs. We review this construction in detail, flanked by a range of examples. Furthermore we explain how open/closed TQFTs are equivalent to Calabi-Yau categories and the Cardy condition, and how to extract such data from pivotal 2-categories.

# Contents

# 1 Introduction

Defects in field theories describe various interesting phenomena in physics, and their conceptual significance is becoming increasingly apparent in mathematics. In particular, it is just as natural and useful to consider defects in topological quantum field theories as it is to consider morphisms in categories.

In physics, topological defects are lower-dimensional regions in spacetime where something special is going on. They are 'defective' in the sense that they are different in nature and/or substance from their surroundings. They are 'topological' because geometric details like metrics are not necessary to characterise them; typically this goes hand in hand with a high degree of stability. This loose description of topological defects captures the common denominator of various phenomena in cosmology (cosmic strings), fluid dynamics (hydrodynamic solitons), condensed matter (domain walls in ferromagnets), protein folding (topological frustration), and topological quantum computing.

In mathematical physics and pure mathematics, topological defects take on a more conceptual role. They relate and compare different quantum field theories, and they provide a unifying perspective. In particular, symmetries of a given QFT and dualities between distinct theories (such as mirror symmetry) are special examples of topological defects.

The purpose of these notes is to give an introduction to topological defects in 2-dimensional topological quantum field theory. Such 'defect TQFTs' are defined in the spirit of Atiyah and Segal, as functors on certain decorated bordism categories. Recall that a 2-dimensional *closed* TQFT is a symmetric monoidal functor $\mathcal{Z}^{\mathrm{c}} : \mathrm{Bord}_2 \to \mathrm{Vect}_{\Bbbk}$, where $\mathrm{Bord}_2$ has oriented circles $S^1$ as objects and oriented bordism classes as morphisms.[1] As we will discuss in detail in the next section, a *defect* TQFT is obtained by enlarging the bordism category: both objects and morphisms may have submanifolds of codimension 1, decorated by certain data. An example of such a 'defect bordism' is the decorated surface

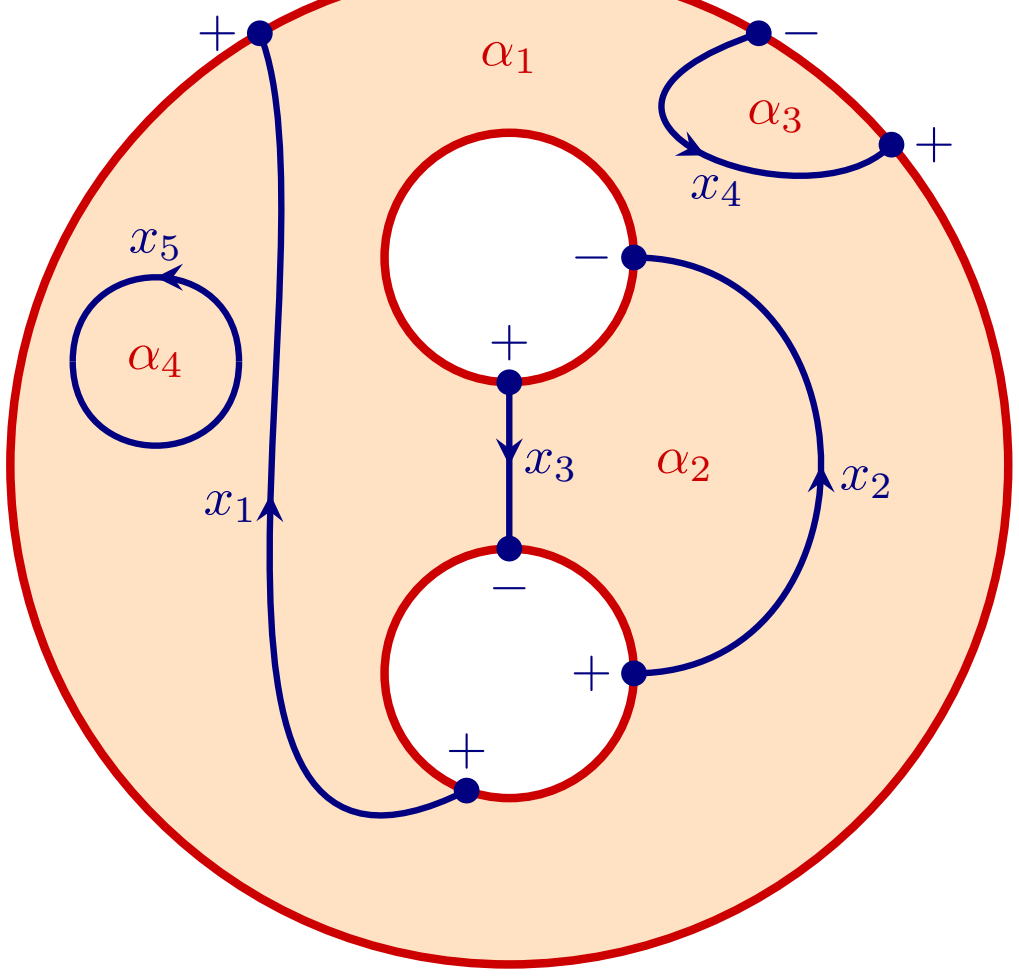

(1.1)

Closed TQFTs are recovered as special cases of defect TQFTs by restricting to trivial decorations only. For instance, if we forget all decorations in (1.1) we obtain the familiar pair-of-pants

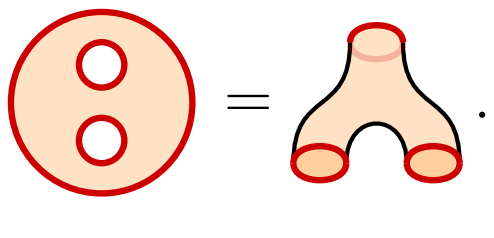

More generally, *open/closed* TQFTs $\mathcal{Z}^{\mathrm{oc}}$ can also be subsumed under the umbrella of defect TQFTs. As we will review independently in Section 3.2, by definition $\mathcal{Z}^{\mathrm{oc}}$ does not only act on $S^1$, but also on intervals $I_{ab}$ with endpoints labelled by elements $a, b$ of some chosen set $B$ of 'boundary conditions'. We will also see in which precise way open/closed TQFTs are the special cases of defect TQFTs whose only defect lines are boundaries. Hence defect TQFTs reside on top of the sequence

$$\text{closed TQFTs} \subsetneq \text{open/closed TQFTs} \subsetneq \text{defect TQFTs}. \tag{1.2}$$

[1] Throughout these notes we assume that every manifold comes with an orientation.

It is well-known that 2-dimensional closed TQFTs $\mathcal{Z}^{\mathrm{c}}$ are equivalently described by commutative Frobenius algebras, namely the vector space $\mathcal{Z}^{\mathrm{c}}(S^1)$ with (co)multiplication coming from the (co)pair-of-pants. Similarly, open/closed TQFTs $\mathcal{Z}^{\mathrm{oc}}$ naturally give rise to a category whose objects are the boundary conditions $a, b, \ldots \in B$, and whose morphism spaces are what $\mathcal{Z}^{\mathrm{oc}}$ assigns to the decorated intervals $I_{ab}$, cf. Section 3.2. As a generalisation of these facts we will explain how every 2-dimensional defect TQFT $\mathcal{Z}$ naturally gives rise to a certain 2-category $\mathcal{B}_{\mathcal{Z}}$ [DKR]. Its objects are to be thought of as closed TQFTs, its 1-morphisms correspond to line defects, and its 2-morphisms are operators associated to intersection points of defect lines. Thus in a nutshell:

$$\begin{aligned} \text{closed TQFT} &\Longrightarrow \text{algebra} \\ \text{open/closed TQFT} &\Longrightarrow \text{category} \\ \text{defect TQFT} &\Longrightarrow \text{2-category} \end{aligned}$$

One motivation to consider the 2-category $\mathcal{B}_{\mathcal{Z}}$ is the following simple motto:

*Correlators* of the theory $\mathcal{Z}$ are *string diagrams* in the 2-category $\mathcal{B}_{\mathcal{Z}}$.

Indeed, we will explain how topological correlators rigorously translate into the algebraic setting of $\mathcal{B}_{\mathcal{Z}}$, where they can be easily manipulated and computed. As specific examples we will discuss sphere and disc correlators in detail.

Another reason to adopt a higher-categorical language to study defect TQFTs is that the conceptual clarity and bird's eye view make known symmetries and dualities appear more natural. Even better, the algebraic framework of $\mathcal{B}_{\mathcal{Z}}$ paves the way to new dualities and new equivalences of categories [FFRS, CR, CRCR, CQV].

The remainder of these notes is organised as follows. In Section 2.1 we introduce the 2-dimensional defect bordism category $\mathrm{Bord}_2^{\mathrm{def}}(\mathbb{D})$. Then we define TQFTs as symmetric monoidal functors $\mathcal{Z} : \mathrm{Bord}_2^{\mathrm{def}}(\mathbb{D}) \to \mathrm{Vect}_{\Bbbk}$ and explain how closed and open/closed TQFTs are special cases. (An independent introduction to open/closed TQFT is given in Section 3.2.) In Section 2.2 we review the basic notions for 2-categories which are needed for the construction of a 'pivotal' 2-category $\mathcal{B}_{\mathcal{Z}}$ for every defect TQFT $\mathcal{Z}$ in Section 2.3. Then in Section 2.4 we discuss various examples.

Section 3 takes pivotal 2-categories as the starting point and shows how to construct open/closed TQFTs from them, under two natural assumptions. The purely closed case is discussed in Section 3.1, where for every object in a 2-category satisfying Assumption 3.1 we construct a commutative Frobenius algebra. In the parenthetic Section 3.2 we review open/closed TQFTs and discuss how they are algebraically encoded in terms of Calabi-Yau categories and the Cardy condition. Finally in Section 3.3 we extract such algebraic data from every object in a 2-category satisfying Assumption 3.1 and 3.7. Thus we make sense of the purely algebraic version of (1.2):

$$\text{Frobenius algebras} \subsetneq \text{Calabi-Yau categories} \subsetneq \text{pivotal 2-categories.}$$

### Acknowledgements

I thank F. Montiel Montoya, D. Murfet, P. Pandit, A. G. Passegger, D. Plencner, G. Schaumann, D. Scherl, D. Tubbenhauer, P. Wedrich, and K. Wehrheim for comments and discussions. This work is partially supported by a grant from the Simons Foundation, and from the Austrian Science Fund (FWF): P 27513-N27.

# 2 Defect TQFT

## 2.1 Functorial definition

To describe an open/closed TQFT one has to specify the set of boundary conditions. Similarly, to describe an $n$-dimensional defect TQFT we need sets $D_j$ whose elements label the $j$-dimensional defects in bordisms for $j \in \{1, 2, \dots, n\}$. Note that for *topological* theories we do not require a set $D_0$ as input data to label points; we will see how the defect TQFT itself computes the set $D_0$.

In addition to the defect label sets $D_j$, we also need a set of maps $\mathcal{D}$ between them to encode how defects of different dimensions are allowed to meet – e. g. which labels in $D_{n-1}$ may occur on domain walls between two given $n$-dimensional regions, or what the neighbourhood around a $D_1$-labelled defect line can look like. We collectively refer to the sets $D_j$ and maps $\mathcal{D}$ as *defect data* $\mathbb{D}$. We will momentarily describe $\mathbb{D}$ in detail for the case $n = 2$.

An $n$-dimensional defect TQFT naturally gives rise to an $(n+1)$-layered structure: The basic layer ('objects') is formed by the set $D_n$. The next layer ('1-morphisms') mediates between elements of $D_n$ (via the maps $\mathcal{D}$), and is comprised of lists of elements of $D_{n-1}$. Then there is a layer ('2-morphisms') made of elements in $D_{n-2}$, which is between 1-morphisms. This goes on until we have $(n-1)$-morphisms built from $D_1$, and finally there are $n$-morphisms which are computed from the TQFT itself. Unsurprisingly, this $(n + 1)$-layered structure turns out to be an '$n$-category'. The case $n = 2$ was first worked out in [DKR], and we will now review it in detail.[2]

[2] The case $n = 3$ was worked out in [BMS]. The details for $n > 3$ have not been worked out.

We assume the reader is familiar with Atiyah's definition of 2-dimensional closed TQFTs as symmetric monoidal functors $\mathrm{Bord}_2 \to \mathrm{Vect}_{\Bbbk}$ as reviewed in [Koc]. Furthermore, we will exclusively consider oriented TQFTs; hence all the circles and bordisms below implicitly come with an orientation. A 2-dimensional defect TQFT is a generalisation where the bordism category $\mathrm{Bord}_2$ is enlarged: by definition a *defect TQFT* is a symmetric monoidal functor

$$\mathcal{Z} : \mathrm{Bord}_2^{\mathrm{def}}(D_1, D_2, s, t) \longrightarrow \mathrm{Vect}_{\Bbbk} \, .$$

What is $\mathrm{Bord}_2^{\mathrm{def}}(D_1, D_2, s, t)$? First of all, $D_1$ and $D_2$ are any two chosen sets, which will label 1-dimensional lines and 2-dimensional regions on bordisms, respectively. The *source* and *target maps*

$$s, t : D_1 \longrightarrow D_2$$

tell us how $D_2$-labelled regions may meet at $D_1$-labelled lines. To wit, in the defect bordisms introduced below, the region to the right of a line labelled by $x \in D_1$ must be labelled by $s(x) \in D_2$, and the label on the left is $t(x) \in D_2$:

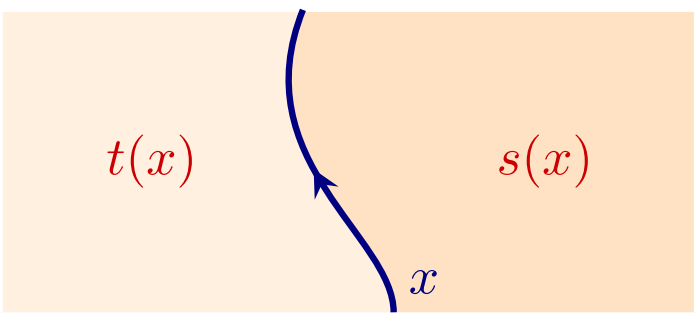

We now use the chosen *defect data*

$$\mathbb{D} := (D_1, D_2, s, t)$$

to decorate the objects and morphisms in $\mathrm{Bord}_2^{\mathrm{def}}(\mathbb{D})$. An *object* in $\mathrm{Bord}_2^{\mathrm{def}}(\mathbb{D})$ is a disjoint union of circles $S^1$ with finitely many points $p \in S^1 \setminus \{-1\}$ labelled by pairs $(x, \varepsilon)$ with $x \in D_1$ and $\varepsilon \in \{\pm\}$, and line segments between such points $p$ are labelled by elements in $D_2$. Such decorated circles are called *defect circles*. An example of an object in $\mathrm{Bord}_2^{\mathrm{def}}(\mathbb{D})$ is the disjoint union

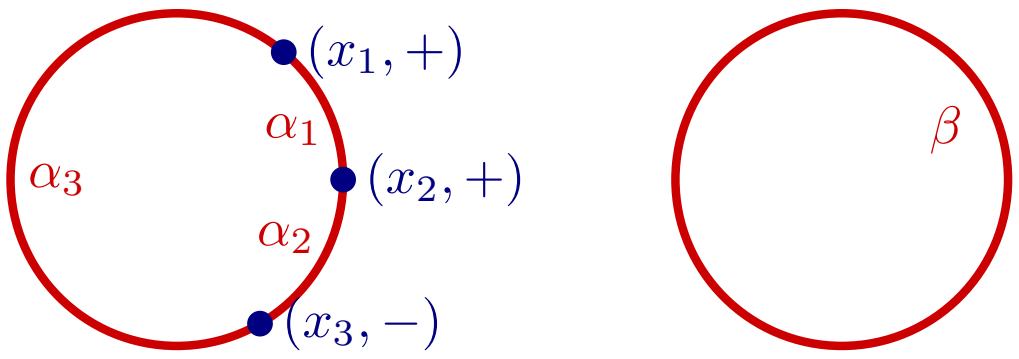

where the second defect circle has no marked point and is decorated with $\beta \in D_2$.

A *morphism* in $\mathrm{Bord}_2^{\mathrm{def}}(\mathbb{D})$ is either a permutation of the labels on a given defect circle, or a *defect bordism* class. A defect bordism is an ordinary bordism $\Sigma$ in $\mathrm{Bord}_2$ together with an oriented 1-dimensional submanifold $\Sigma_1$. This submanifold may have nonempty boundary $\partial\Sigma_1$, but we require $\partial\Sigma_1$ to lie in the boundary of the bordism: $\partial\Sigma_1 \subset \partial\Sigma$. Each connected component ('defect line') of $\Sigma_1$ is labelled by an element in $D_1$, while the components ('phases') of $\Sigma \setminus \Sigma_1$ are labelled by elements in $D_2$ such that the phase to the right (respectively left) of an $x$-labelled defect line is decorated by $s(x)$ (respectively $t(x)$). Defect lines may meet the boundary of $\Sigma$ only transversally, and only at marked points $p$ of the associated defect circle. If $p$ is decorated by $(x, \varepsilon)$ and sits on the in-going boundary, then the defect line touching $p$ must also be labelled by $x$, and it must be oriented away from (respectively towards) the boundary if $\varepsilon = +$ (respectively $\varepsilon = -$). If $p$ sits on the out-going boundary then the role of the sign $\varepsilon$ is reversed. Finally, the $D_2$-labels of line segments in objects in $\mathrm{Bord}_2^{\mathrm{def}}(\mathbb{D})$ must coincide with those of their adjacent phases in $\Sigma \setminus \Sigma_1$.

Two defect bordisms belong to the same class if there exists an isotopy between them whose restriction to defect lines is a bijection, and if the defect labels are the same. The condition that marked points on defect circles cannot sit at $-1 \in S^1$ ensures that rotating all marked points by $2\pi$ is disallowed and hence cannot give rise to non-trivial endomorphisms.

An example of a defect bordism is the decorated surface

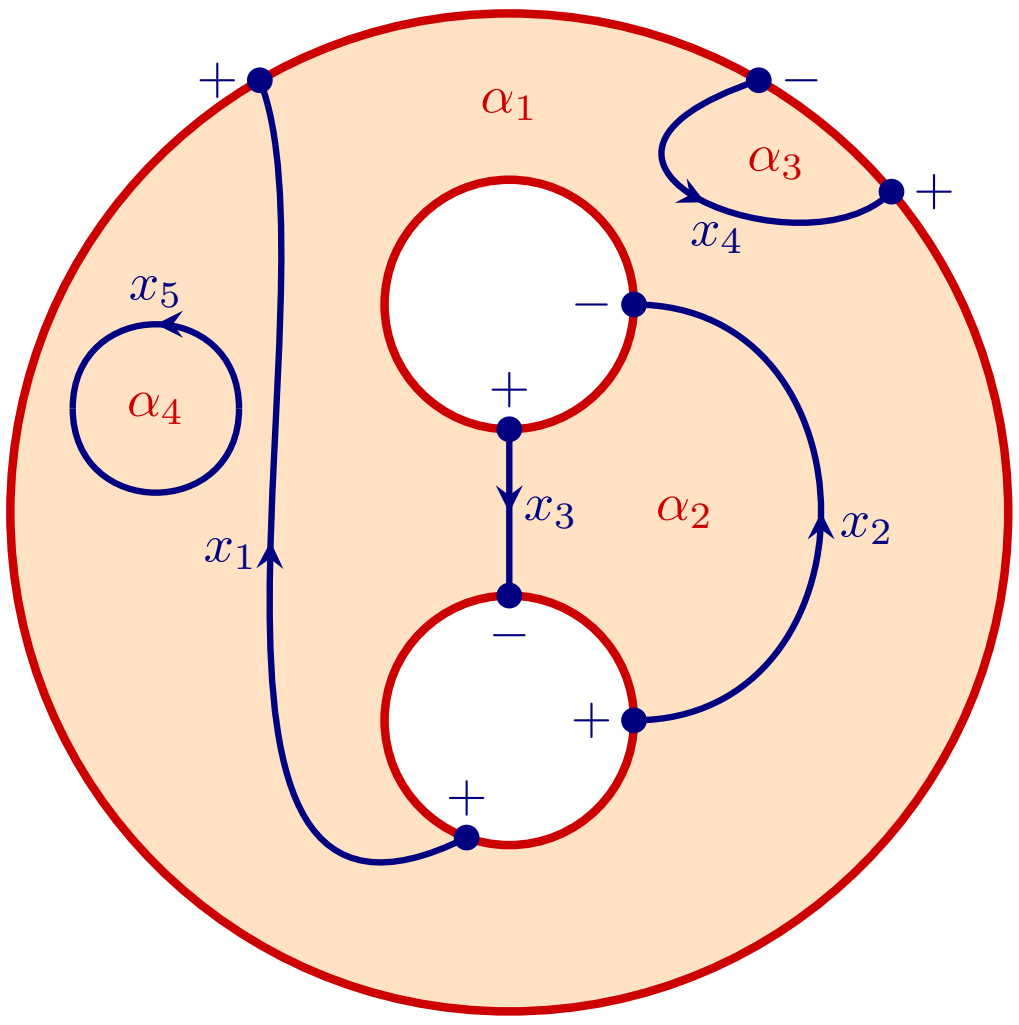

(2.1)

where we choose to view the two inner circles as in-coming boundaries, and the outer circle as out-going. Hence this 'pair-of-pants' represents a morphism

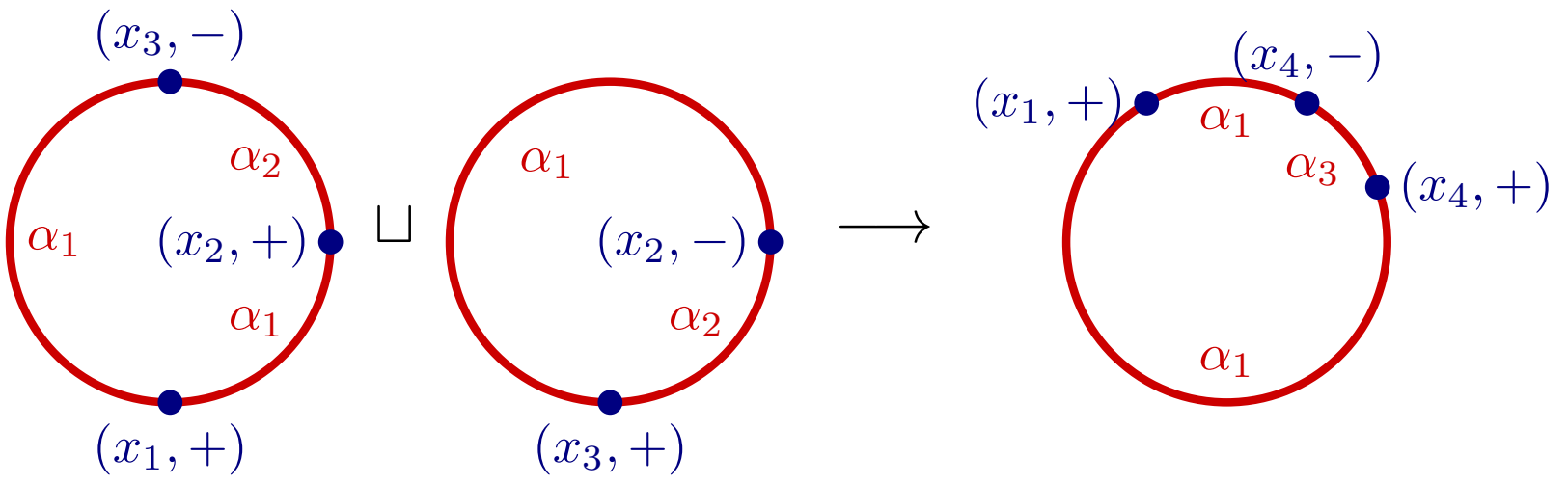

.

As in (2.1), we sometimes do not decorate boundaries with $D_1, D_2$ when drawing bordisms, as these decorations can be inferred from those of the interior.

### 2.1.1 Open/closed TQFTs as defect TQFTs

A closed TQFT is the special case of a 'defect TQFT without defects'. More precisely, if we choose $D_1 = \emptyset$ ('no defect lines') and $D_2 = \{\bullet\}$ ('only a single phase') then $s, t$ are trivial, and forgetting the label $\bullet$ attached to every bordism is an equivalence $\mathrm{Bord}_2^{\mathrm{def}}(\emptyset, \{\bullet\}) \to \mathrm{Bord}_2$.

Another special case are open/closed TQFTs $\mathrm{Bord}_2^{\mathrm{oc}}(B) \to \mathrm{Vect}_{\Bbbk}$, where objects in the bordism category may also involve intervals with endpoints labelled by elements of a set (of 'boundary conditions') $B$, cf. Section 3.2. Indeed, let us choose

$$D_1 = B \, , \quad D_2 = \{\bullet, \circ\} \quad \text{with} \quad s(D_1) = \{\circ\} \, , \quad t(D_1) = \{\bullet\} \, ,$$

where we interpret $\circ$ as the label for the 'trivial theory' (whose associated Frobenius algebra is $\Bbbk$). Then there is a forgetful functor $U : \mathrm{Bord}_2^{\mathrm{def}}(B, \{\bullet, \circ\}, s, t) \to \mathrm{Bord}_2^{\mathrm{oc}}(B)$ which erases all line segments and phases labelled by $\circ$, and which forgets the label $\bullet$ for the remaining line segments and phases. For example, a flat pair-of-pants in $\mathrm{Bord}_2^{\mathrm{oc}}(B)$ is obtained as

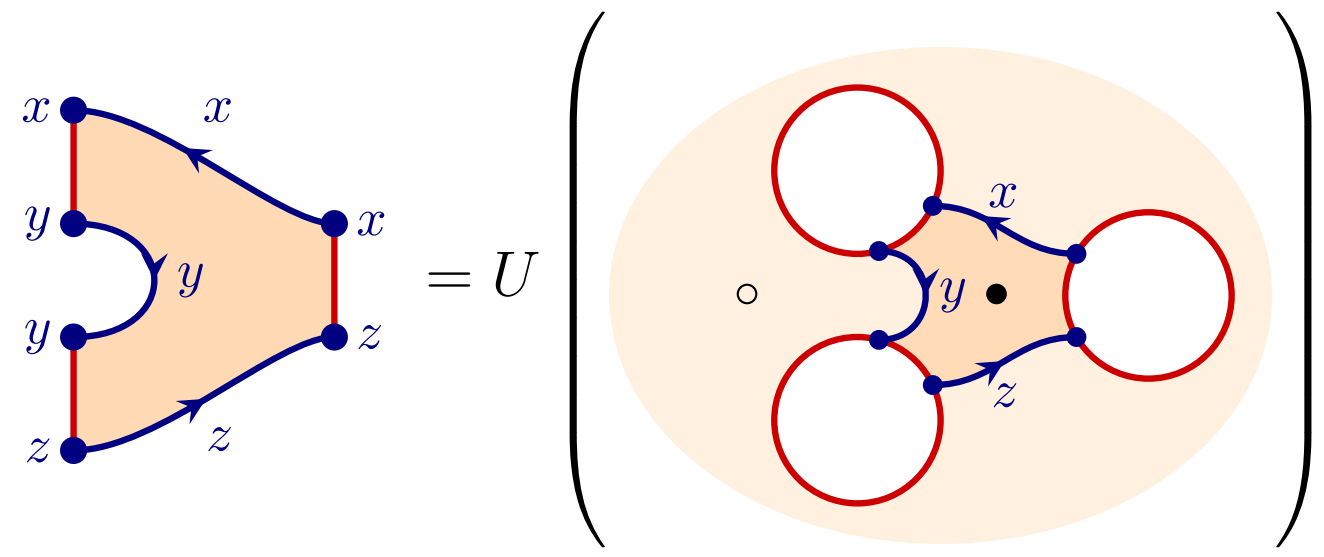

where the $\circ$-labelled phase may close to any bordism, e. g. the sphere.

A 2-dimensional closed TQFT $\mathcal{Z}^{\mathrm{c}}$ is equivalently described by a commutative Frobenius algebra [Koc]. Similarly, an open/closed TQFT $\mathcal{Z}^{\mathrm{oc}}$ is equivalent to a commutative Frobenius algebra (what $\mathcal{Z}^{\mathrm{oc}}$ assigns to a circle) and a Calabi-Yau category (whose set of objects is $B$, and whose Hom spaces are what $\mathcal{Z}^{\mathrm{oc}}$ assigns to intervals), together with certain maps [Laz, AN, MS, LP1]. Below in Section 3.2 we will review this in more detail, and in Section 3.3 we will explain how all these algebraic structures are naturally obtained from the general perspective of defect TQFT. However, we first need to bring the process

$$\overset{\text{closed}}{\text{TQFT}} \approx \text{Frobenius algebra} \xrightarrow[\text{structure}]{\text{add boundary}} \overset{\text{open/closed}}{\text{TQFT}} \approx \text{Calabi-Yau category}$$

to its logical conclusion and show how the expectation

$$\overset{\text{defect}}{\text{TQFT}} \approx \text{pivotal 2-category}$$

can be made rigorous.

## 2.2 Pivotal 2-categories

We start with a brief review of 2-categories and their graphical calculus; for a more detailed account we refer to [Ben, KS, Lau2]. Let us first recall that every $\Bbbk$-algebra $A$ can be viewed as the endomorphism space $\mathrm{End}(*)$ of a $\Bbbk$-linear category with a single object $*$. In this sense a category generalises the idea of 'many algebras together'.

A 2-category is 'many monoidal categories together'. The precise definition is that a *2-category* is a category enriched over the category of small categories. This means that for any two objects $\alpha, \beta$ of a 2-category $\mathcal{B}$, there is a category $\mathcal{B}(\alpha, \beta)$ whose objects are called *1-morphisms* from $\alpha$ to $\beta$, and whose morphisms are called *2-morphisms*. The composition of two 2-morphisms in $\mathcal{B}(\alpha, \beta)$ is called *vertical composition*. Since 1-morphisms can also be composed (being morphisms in $\mathcal{B}$) there are functors

$$\otimes : \mathcal{B}(\beta, \gamma) \times \mathcal{B}(\alpha, \beta) \longrightarrow \mathcal{B}(\alpha, \gamma)$$

which we refer to as *horizontal composition*. A good example of a 2-category is that whose objects, 1- and 2-morphisms are small categories, functors and natural transformations, respectively.

The attributes 'vertical' and 'horizontal' for the two types of composition in a 2-category $\mathcal{B}$ derive from the *graphical calculus*, which allows to perform computations in $\mathcal{B}$ in terms of so-called *string diagrams*, according to the following rules:

- Objects in $\mathcal{B}$ label 2-dimensional regions in the plane.
- 1-morphisms $X : \alpha \to \beta$ label smooth lines with $\alpha$ to the right and $\beta$ to the left. The lines must be progressive, i. e. at no point on the line may the tangent vector have zero component in upward direction, cf. [BMS, Def. 2.8]. For example:

  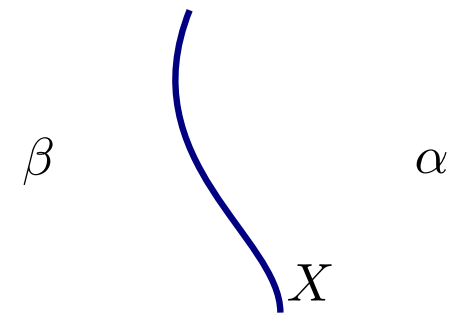
  

  Such an $X$-line is also identified with the unit 2-morphism $1_X \in \mathrm{End}(X)$.

- 2-morphisms $\Phi : X \to Y$ label vertices on a line with label $X$ below and label $Y$ above the vertex, respectively:

  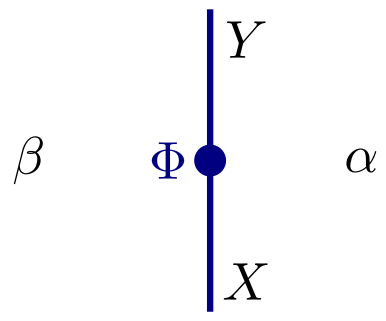
  
  .

  This diagram is identified with $\Phi$.

- Vertical composition of $\Phi : X \to Y$ and $\Psi : Y \to Z$ really is vertical, read from bottom to top:

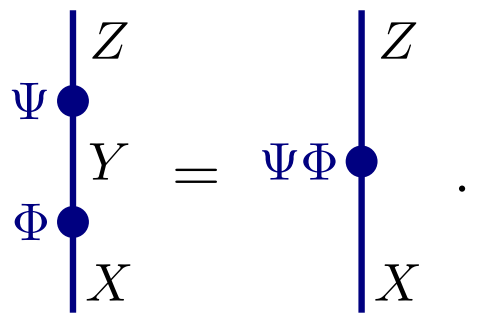

  Note that as above we sometimes suppress labels for 2-dimensional regions.

- Horizontal composition really is horizontal, read from right to left:

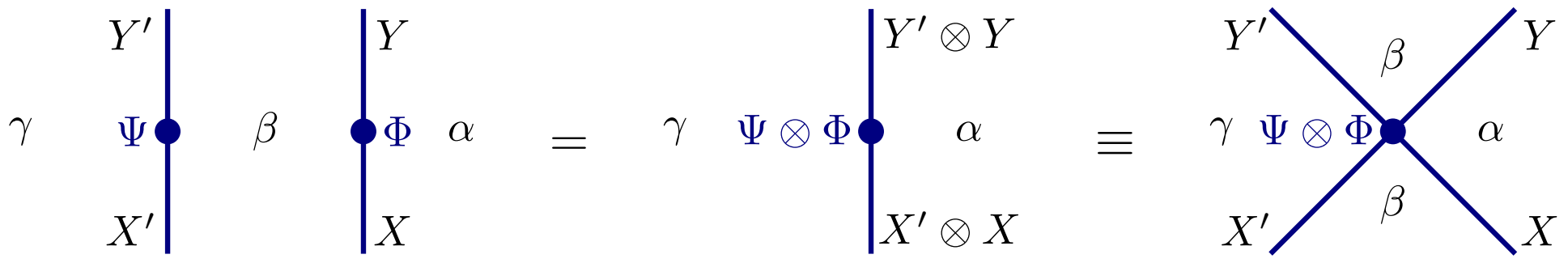

  Here we used the additional rule that elements in the set of 2-morphisms $\mathrm{Hom}(X_1 \otimes \cdots \otimes X_m, Y_1 \otimes \cdots \otimes Y_n)$ may be depicted as vertices with $m$ in-coming $X_i$-lines and $n$ out-going $Y_j$-lines. Consistency with $X = X \otimes 1_\alpha$ then demands that the unit 1-morphisms $1_\alpha$ can be represented by invisible lines.

Hence every progressive string diagram represents a 2-morphism by reading it from bottom to top and from right to left. For example

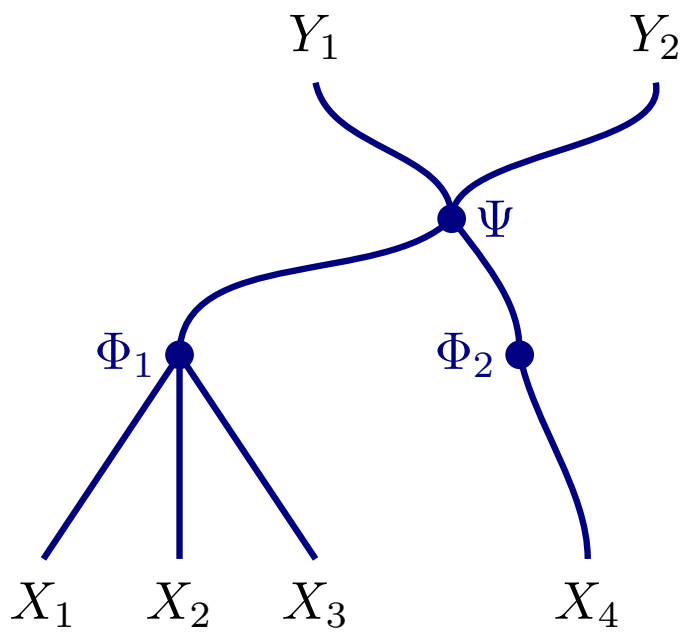

$$= \Psi \circ (\Phi_1 \otimes \Phi_2) : X_1 \otimes X_2 \otimes X_3 \otimes X_4 \longrightarrow Y_1 \otimes Y_2 \,. \quad (2.2)$$

But what if the loci of the lines or vertices vary a little? Does

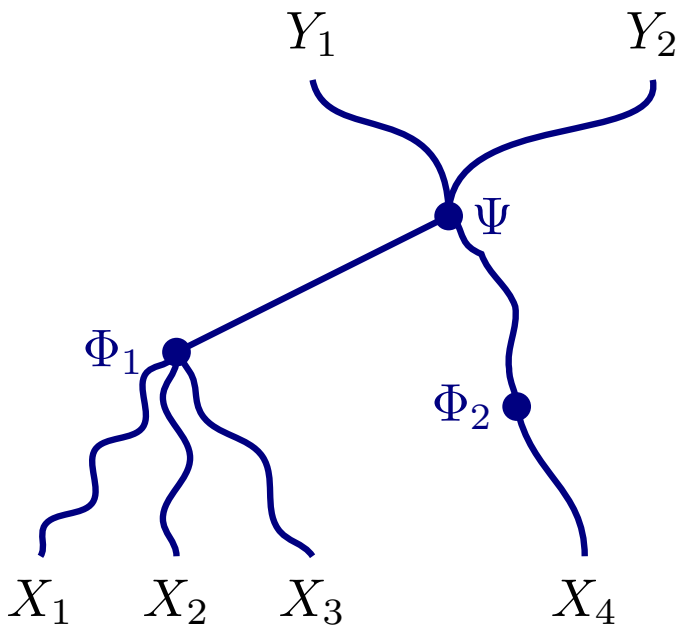

represent the same 2-morphism as (2.2)? It better should, and this is guaranteed by the following rule:

- String diagrams which are related by progressive isotopies represent the same 2-morphism [BMS, Sect. 2.2].

This in particular implies the interchange law

$(1 \otimes \Phi) \circ (\Psi \otimes 1) =$ 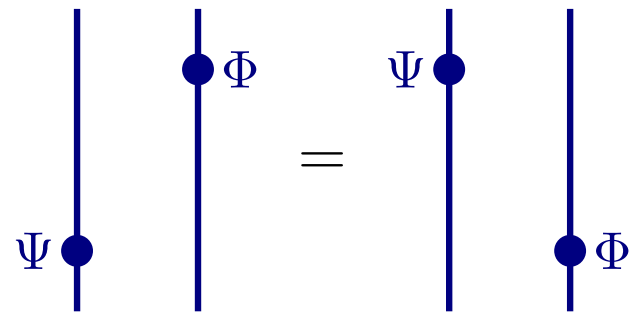 $= (\Psi \otimes 1) \circ (1 \otimes \Phi)$ (2.3)

which is a consequence of the functoriality of $\otimes$.

String diagrams for 2-categories are already reminiscent of local patches on defect bordisms. To make this relation precise we need to enlarge the type of string diagrams we consider, by giving an orientation to every line, and by allowing them to make 'U-turns'. In order to continue to have a 1-to-1 relation between isotopy classes of string diagrams and 2-morphisms, we need however to consider 2-categories with additional structure: 2-categories 'with adjoints'.

Given a 1-morphism $X : \alpha \to \beta$ in a 2-category $\mathcal{B}$, we say that $X$ has a *left adjoint* if there is ${}^\dagger X \in \mathcal{B}(\beta, \alpha)$ together with 2-morphisms

$$\mathrm{ev}_X : {}^\dagger X \otimes X \longrightarrow 1_\alpha \, , \quad \mathrm{coev}_X : 1_\beta :\longrightarrow X \otimes {}^\dagger X$$

subject to the conditions

$$\big(1_X \otimes \mathrm{ev}_X\big) \circ \big(\mathrm{coev}_X \otimes 1_X\big) = 1_X \, , \quad \big(\mathrm{ev}_X \otimes 1_{{}^\dagger X}\big) \circ \big(1_{{}^\dagger X} \otimes \mathrm{coev}_X\big) = 1_{{}^\dagger X} \, . \tag{2.4}$$

Since the *adjunction maps* $\mathrm{ev}_X, \mathrm{coev}_X$ are something special, they deserve special diagrammatic notation:

$\mathrm{ev}_X =$ 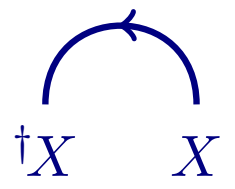 , $\mathrm{coev}_X =$ 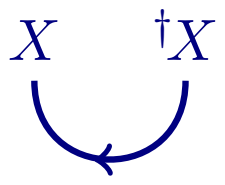 . (2.5)

Here we are using our final rule for string diagrams:

- Lines for objects $X$ with an adjoint come with an orientation. Such lines need not be progressive, but the only non-progressive parts must be diagrams for adjunction maps as in (2.5) or (2.6). Upward-oriented line segments are labelled $X$, downward-oriented segments are labelled ${}^\dagger X$ or $X^\dagger$.

In diagrammatic language the conditions (2.4) are easy to remember: they are called *Zorro moves* and state that 'lines may be straightened out':

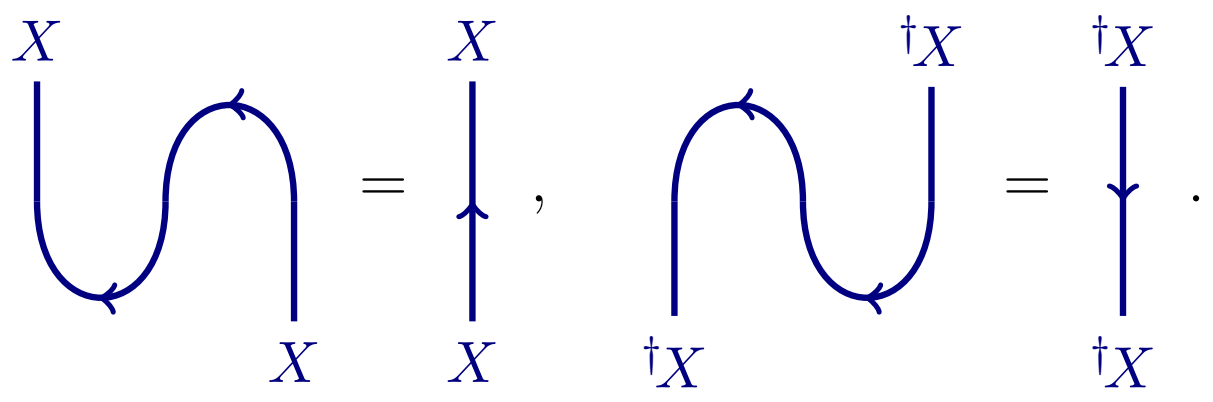

Similarly, the *right adjoint* of $X \in \mathcal{B}(\alpha, \beta)$ is $X^\dagger \in \mathcal{B}(\beta, \alpha)$ together with adjunction maps

$\widetilde{\mathrm{ev}}_X =$ 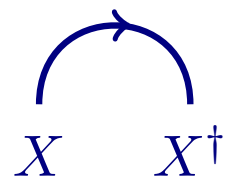

$: X \otimes X^\dagger \longrightarrow 1_\beta \,, \quad \widetilde{\mathrm{coev}}_X =$ 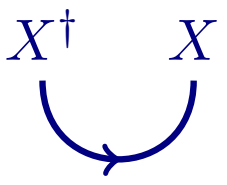

$: 1_\alpha \longrightarrow X^\dagger \otimes X$ (2.6)

that satisfy the Zorro moves

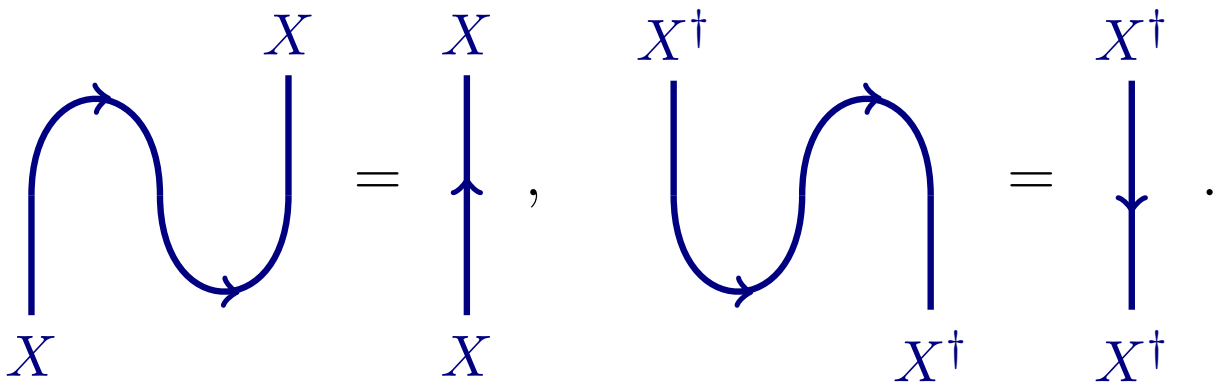

If every 1-morphism in $\mathcal{B}$ has a left and a right adjoint, we say that $\mathcal{B}$ *has adjoints*.

While it is not difficult to prove that left and right adjoints are unique up to isomorphism, ${}^\dagger X$ need not be isomorphic to $X^\dagger$ in general. However, since from the perspective of TQFT taking the adjoint corresponds to orientation reversal on defect lines, we expect ${}^\dagger X = X^\dagger$ in the 2-categories we will construct from defect TQFTs. More precisely, we will encounter *pivotal* 2-categories, which by definition have adjoints with ${}^\dagger X = X^\dagger$ for all 1-morphisms $X$, and where both adjoints for 2-morphisms are identified as well. This means we require the identities

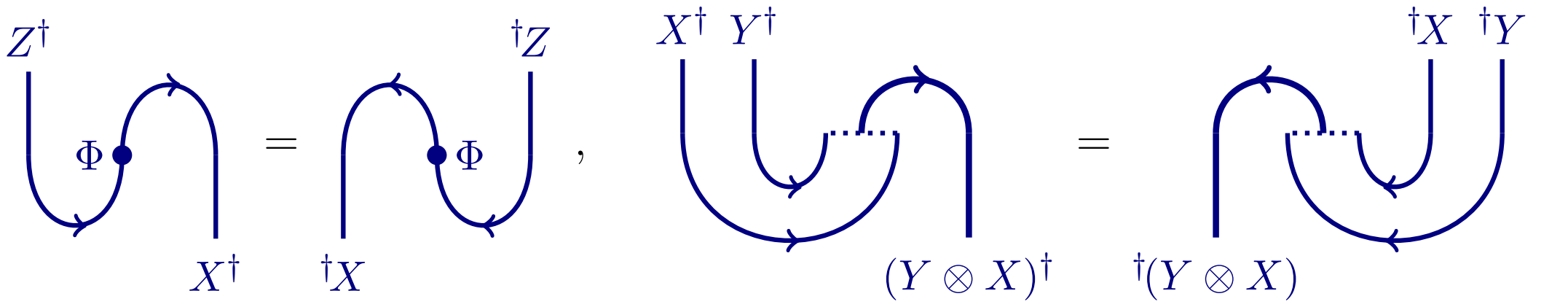

(2.7)

whenever these diagrams make sense.

If the left and right adjoints of a 1-morphism $X \in \mathcal{B}(\alpha, \beta)$ coincide, we can compose $\widetilde{\mathrm{coev}}_X$ with $\mathrm{ev}_X$, and $\mathrm{coev}_X$ with $\widetilde{\mathrm{ev}}_X$. These composites are called the *left* and *right quantum dimensions*:

$\dim_{\mathrm{l}}(X) =$ 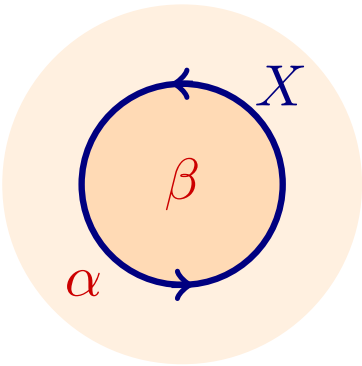

$\in \mathrm{End}(1_\alpha) \,, \quad \dim_{\mathrm{r}}(X) =$ 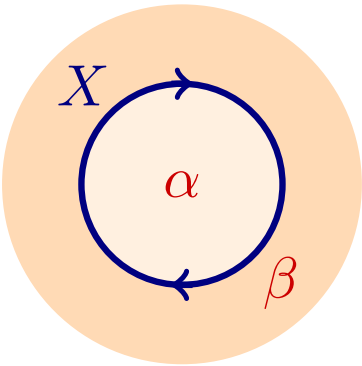

$\in \mathrm{End}(1_\beta) \,.$

(2.8)

More generally, in a pivotal 2-category the *left* and *right traces* of an endomorphism $\Psi \in \mathrm{End}(X)$ are defined to be

$\mathrm{tr_l}(\Psi) =$ 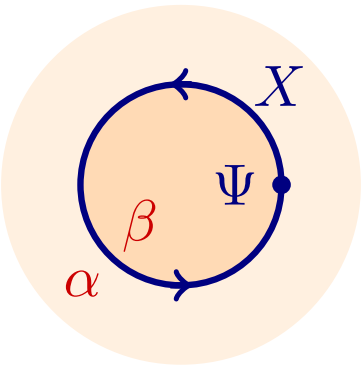

$\in \mathrm{End}(1_\alpha)\,, \quad \mathrm{tr_r}(\Psi) =$ 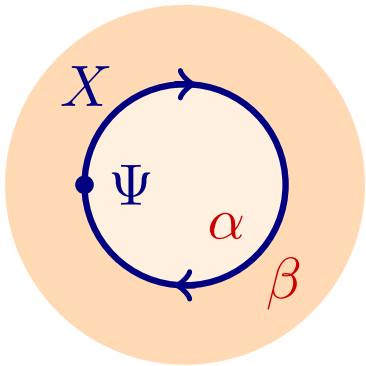

$\in \mathrm{End}(1_\beta)\,. \quad (2.9)$

It follows from the first identity in (2.7) that traces have the expected cyclic property, i. e. $\mathrm{tr_l}(\Phi\Psi) = \mathrm{tr_l}(\Psi\Phi)$ and $\mathrm{tr_r}(\Phi\Psi) = \mathrm{tr_r}(\Psi\Phi)$ for any anti-parallel pair of 2-morphisms $\Phi, \Psi$.

The above notions of adjoints, quantum dimensions and traces generalise the case of finite-dimensional vector spaces $V \in \mathrm{Vect}_\Bbbk$. Indeed, ${}^\dagger V$ and $V^\dagger$ are given by the dual vector space $V^*$, and $\mathrm{ev}_V$ really is the evaluation $V^* \otimes_\Bbbk V \to \Bbbk$, $\varphi \otimes v \mapsto \varphi(v)$. Choosing a basis $\{e_i\}$ of $V$, we set $\mathrm{coev}_V(\lambda) = \lambda \sum_i e_i \otimes e_i^*$ for all $\lambda \in \Bbbk$, and analogously for $\widetilde{\mathrm{ev}}_V$ and $\widetilde{\mathrm{coev}}_V$. Then one easily verifies the Zorro moves, and (2.9) reduces to the ordinary traces of linear operators (under the canonical identification $V^{**} \cong V$). In particular we have $\dim_\mathrm{l}(V) = \dim_\mathrm{r}(V) = \dim_\Bbbk(V)$.

## 2.3 Pivotal 2-categories from defect TQFTs

With the above preparations we can now construct, following [DKR], the pivotal 2-category $\mathcal{B}_\mathcal{Z}$ associated to a defect TQFT

$$\mathcal{Z} : \mathrm{Bord}_2^{\mathrm{def}}(\mathbb{D}) \longrightarrow \mathrm{Vect}_\Bbbk$$

for any set of defect data $\mathbb{D} = (D_1, D_2, s, t)$. It will be convenient to switch between source and target maps depending on orientations, for which we define

$$s(x,+) = s(x)\,, \quad t(x,+) = t(x) \quad \text{and} \quad s(x,-) = t(x)\,, \quad t(x,-) = s(x)$$

for all $x \in D_1$.

Before we start with the construction of $\mathcal{B}_\mathcal{Z}$, let us lead with its interpretation:

$$\begin{aligned} \text{objects} &\;\widehat{=}\; \text{closed TQFTs} \\ \text{1-morphisms} &\;\widehat{=}\; \text{line defects} \\ \text{2-morphisms} &\;\widehat{=}\; \text{`local' operators inserted at defect junctions} \\ \text{vertical composition} &\;\widehat{=}\; \text{operator product} \\ \text{horizontal composition} &\;\widehat{=}\; \text{fusion product} \\ \text{unit 1-morphisms} &\;\widehat{=}\; \text{invisible defects} \\ \text{adjunction} &\;\widehat{=}\; \text{orientation reversal} \end{aligned} \tag{2.10}$$

With this in mind it is no surprise that we define the *objects* of $\mathcal{B}_{\mathcal{Z}}$ to be the label set for 2-dimensional phases on defect bordisms:

$$\mathrm{Obj}(\mathcal{B}_{\mathcal{Z}}) = D_2 \, .$$

Later in Section 3.1 we will extract a commutative Frobenius algebra from $\mathcal{B}_{\mathcal{Z}}$ for every $\alpha \in D_2$, so objects really are closed TQFTs.

The set of *1-morphisms* $\alpha \to \beta$ is defined to be

$$\Big\{ \big((x_1, \varepsilon_1), \dots, (x_n, \varepsilon_n)\big) \in (D_1 \times \{\pm\})^n \,\Big|\, n \geqslant 0 \, , \, s(x_n, \varepsilon_n) = \alpha \, , \, t(x_1, \varepsilon_1) = \beta \, , \\ s(x_i, \varepsilon_i) = t(x_{i+1}, \varepsilon_{i+1}) \text{ for } i \in \{1, 2, \dots, n-1\} \Big\} \, . \tag{2.11}$$

Why? Clearly we want a defect line labelled by $x \in D_1$ to be a 1-morphism between the objects $s(x)$ and $t(x)$:

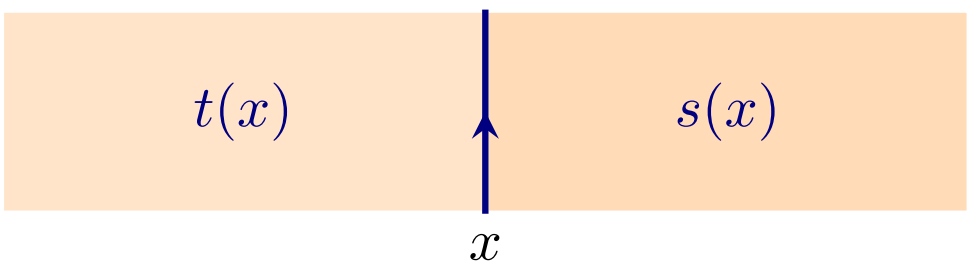

. (2.12)

But if this defect line is 'fused' with another one labelled $y \in D_1$ whose source coincides with the target of $x$,

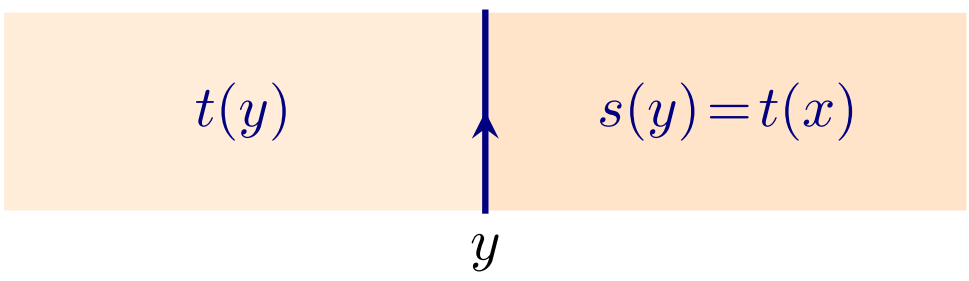

, (2.13)

there is a priori no element '$y \otimes x$' in $D_1$ to label the fusion product of $y$ and $x$. However, we can compose (2.12) and (2.13) to obtain

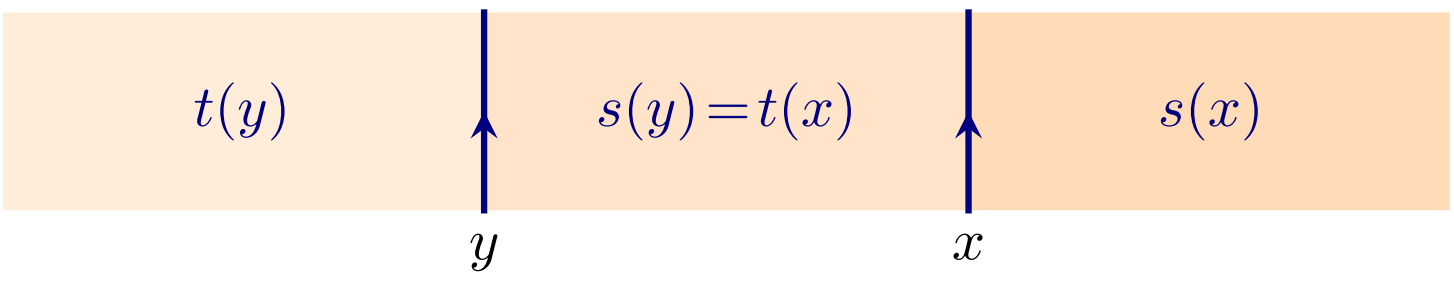

. (2.14)

This picture should be thought of as the composite 1-morphism of $y$ and $x$. And in general, any list of composable defect line labels with orientations $X = ((x_1, \varepsilon_1), \dots, (x_n, \varepsilon_n))$ as in (2.11) is a 1-morphism from $\alpha$ to $\beta$:

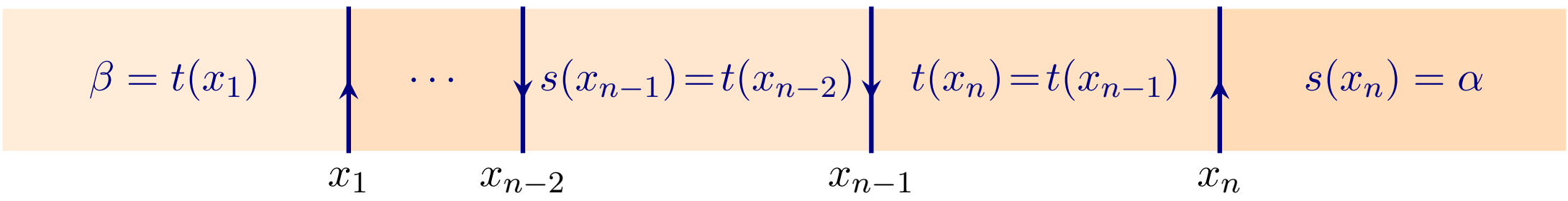

(2.15)

For $X = ((x_1, \varepsilon_1), \dots, (x_n, \varepsilon_n)) : \alpha \to \beta$ and $\widetilde{X} = ((\widetilde{x}_1, \widetilde{\varepsilon}_1), \dots, (\widetilde{x}_m, \widetilde{\varepsilon}_m)) : \beta \to \gamma$ we define *horizontal composition* to be concatenation of lists,

$$\widetilde{X} \otimes X = \big((\widetilde{x}_1, \widetilde{\varepsilon}_1), \dots, (\widetilde{x}_m, \widetilde{\varepsilon}_m), (x_1, \varepsilon_1), \dots, (x_n, \varepsilon_n)\big) : \alpha \longrightarrow \gamma \,.$$

In particular, (2.14) really is the tensor product of (2.12) and (2.13). It is clear that horizontal composition in $\mathcal{B}_{\mathcal{Z}}$ is associative, and the unit 1-morphism $1_\alpha$ is simply the empty sequence ($n = 0$).

The vector space(!) of *2-morphisms* $\mathrm{Hom}(X, Y)$ for $X = ((x_1, \varepsilon_1), \dots, (x_n, \varepsilon_n)) : \alpha \to \beta$ and $Y = ((y_1, \nu_1), \dots, (y_m, \nu_m)) : \alpha \to \beta$ is defined to be

$\mathrm{Hom}(X, Y) = \mathcal{Z}\Bigg($

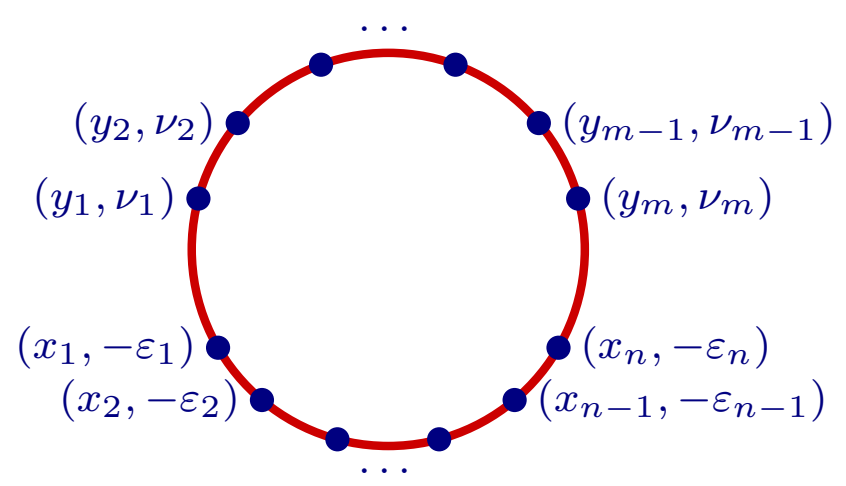

$\Bigg) \,.$ (2.16)

Why? First we note that it is only at this point that we make use of the functor $\mathcal{Z}$. (Objects and 1-morphisms in $\mathcal{B}_{\mathcal{Z}}$ are built only from the defect data $D_1, D_2, s, t$.) Secondly, according to our interpretation (2.10), 2-morphisms should correspond to junction points such as

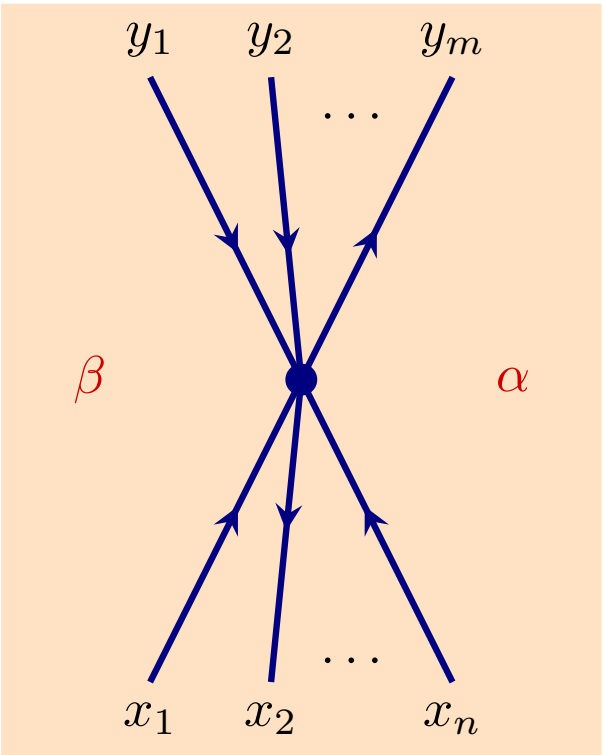

, (2.17)

but we were not provided with a set $D_0$ to label such points. However, we can use $\mathcal{Z}$ to build such a set; in fact it is precisely given by (2.16)! To see this, we cut a little hole around the vertex in (2.17) and keep track of orientations by assigning $+$ to intersection points of lines pointing away from the hole, and $-$ to those pointing in the opposite direction:

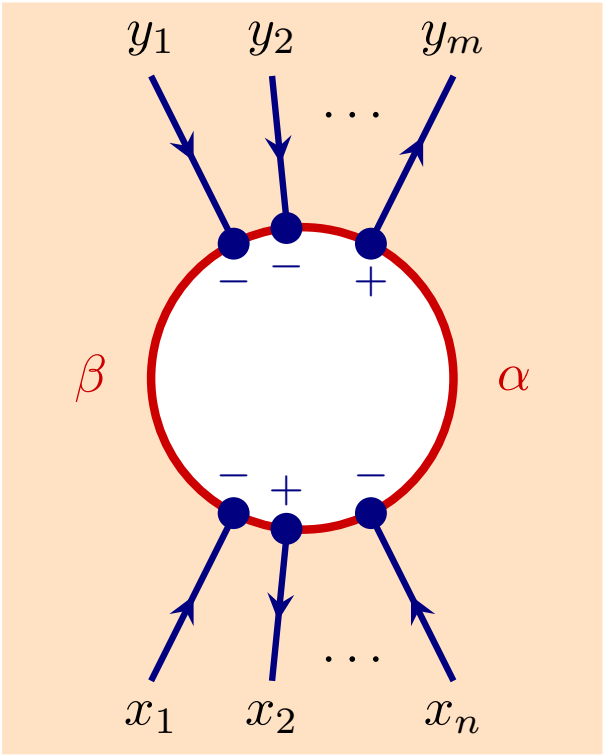

Note that we have not lost any information. But now the boundary of the hole is an object in $\mathrm{Bord}_2^{\mathrm{def}}(\mathbb{D})$. Applying $\mathcal{Z}$ as in (2.16) to it produces the label set for vertices allowed by the TQFT.

Next we define *vertical composition* in $\mathcal{B}_{\mathcal{Z}}$. As in the familiar case of closed TQFT, this product is obtained by applying $\mathcal{Z}$ to a pair-of-pants. In detail, let us consider a 2-morphism $\Phi \in \mathrm{Hom}(X, Y)$ as above, and another 2-morphism $\Psi \in \mathrm{Hom}(Y, Z)$ where $Z = ((z_1, \mu_1), \dots, (z_k, \mu_k)) : \alpha \to \beta$. Then we define

$$\Psi \circ \Phi = \mathcal{Z} \left( \vphantom{\int} \right.$$

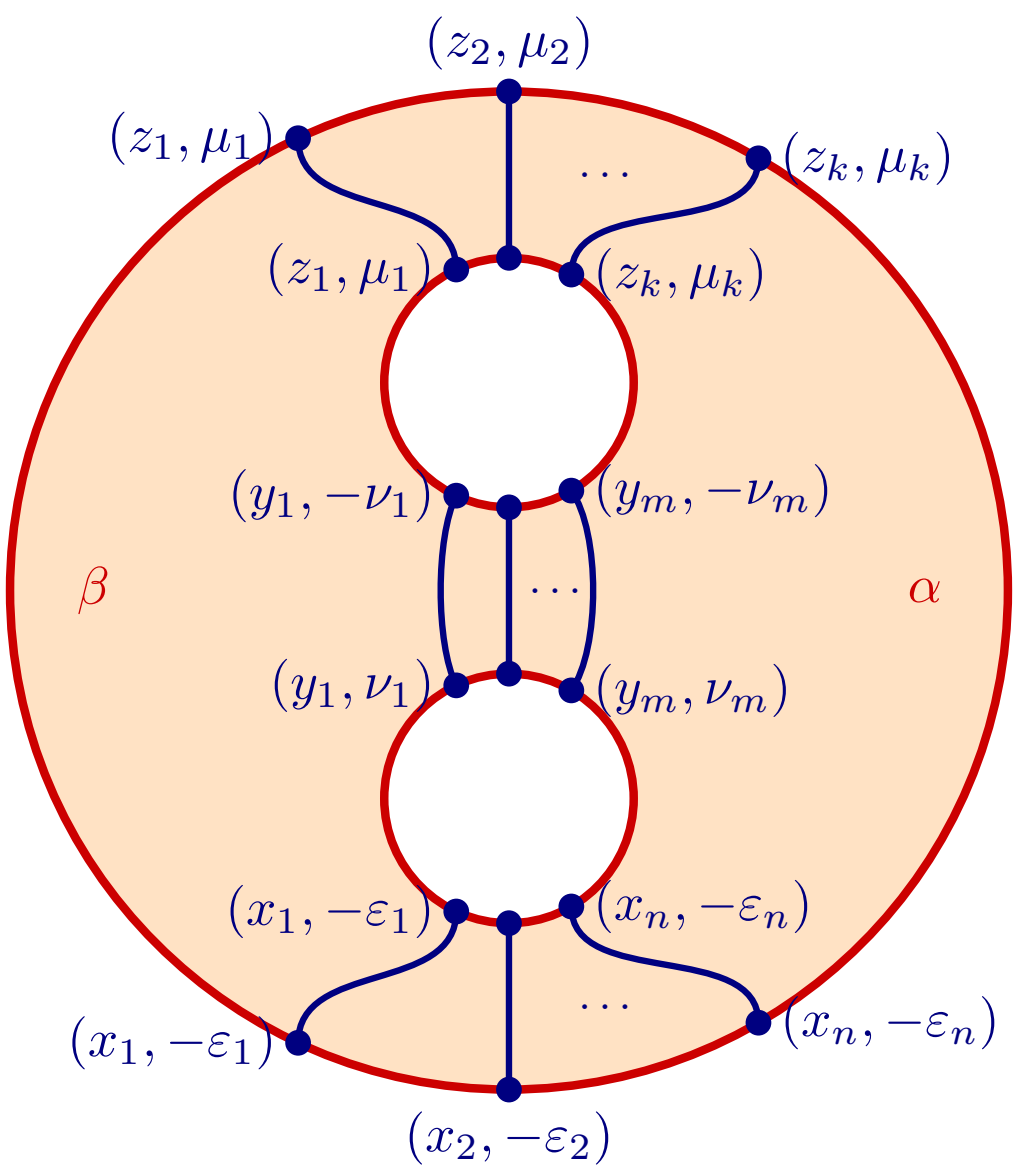

$$\left. \vphantom{\int} \right) \big(\Psi \otimes_{\Bbbk} \Phi\big) \in \mathrm{Hom}(X, Z)$$

where the two inner boundary circles are in-coming, and the outer circle is out-going, and the orientations of defect lines can be inferred from those of their endpoints. Functoriality of $\mathcal{Z}$ implies that this composition is associative, and it is unital with respect to the identity

$$1_X = \mathcal{Z}\left(\right.$$
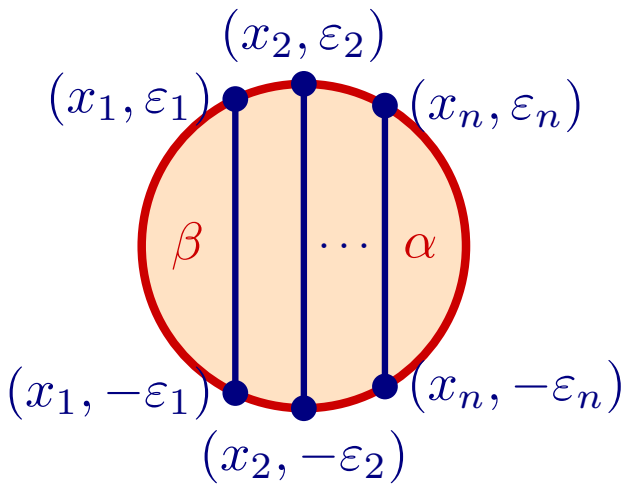

$$\left.\right)(1)\,.$$

Here the disc is viewed as a bordism from the empty set to the boundary circle, so applying $\mathcal{Z}$ gives a linear map $\Bbbk \to \mathrm{End}(X)$.

To complete the construction of $\mathcal{B}_{\mathcal{Z}}$ as a 2-category, we need to define *horizontal composition* of 2-morphisms. Again this involves a pair-of-pants, but this time the defect decoration is different: for $\Phi \in \mathrm{Hom}(X, Y)$ and $\widetilde{\Phi} \in \mathrm{Hom}(\widetilde{X}, \widetilde{Y})$ where $X, Y : \alpha \to \beta$ and $\widetilde{X}, \widetilde{Y} : \beta \to \gamma$, we set

$$\widetilde{\Phi} \otimes \Phi = \mathcal{Z}\left(\right.$$
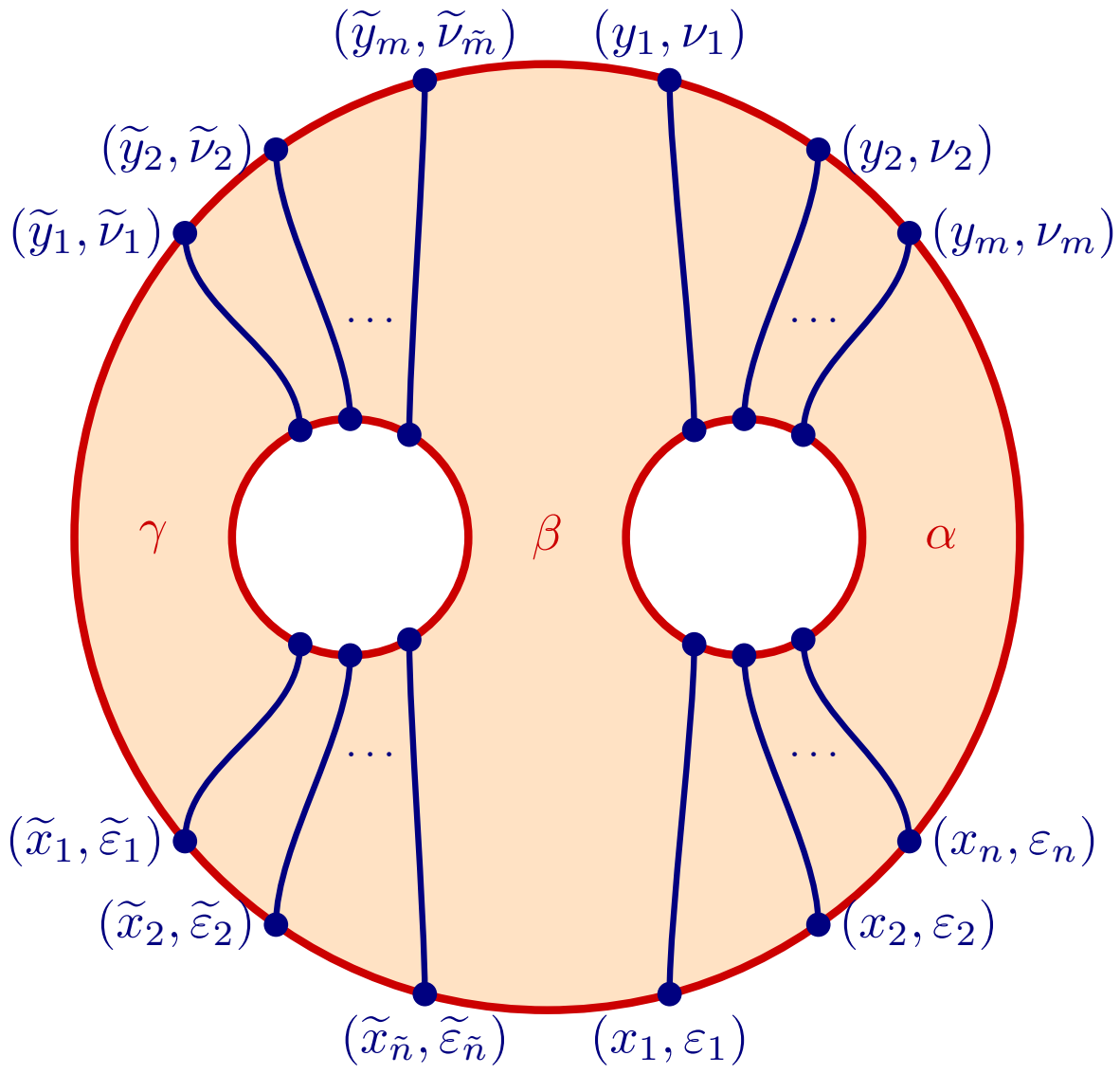

$$\left.\right)\big(\widetilde{\Phi} \otimes_{\Bbbk} \Phi\big)\,.$$

It remains to give $\mathcal{B}_{\mathcal{Z}}$ adjoints and show that it is pivotal. Since adjoints correspond to orientation reversal, adjoints are defined to have opposite signs $\varepsilon$ and reversed order:

$$\begin{aligned} \mathcal{B}_{\mathcal{Z}}(\alpha, \beta) \ni X &= \big((x_1, \varepsilon_1), \dots, (x_n, \varepsilon_n)\big) \\ \Longrightarrow \qquad {}^{\dagger}X \equiv X^{\dagger} &= \big((x_n, -\varepsilon_n), \dots, (x_1, -\varepsilon_1)\big) \in \mathcal{B}_{\mathcal{Z}}(\beta, \alpha)\,. \end{aligned}$$

Note that this was already implicitly used in the definition (2.16), so we have $\mathrm{Hom}(X, Y) = \mathrm{Hom}(1_\beta, Y \otimes X^{\dagger})$.

To exhibit $^\dagger X$ as the left adjoint to $X$ we define the adjunction maps

$\mathrm{ev}_X =$ 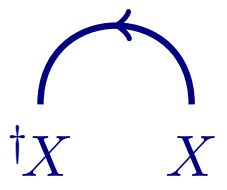

$= \mathcal{Z}\Bigg($ 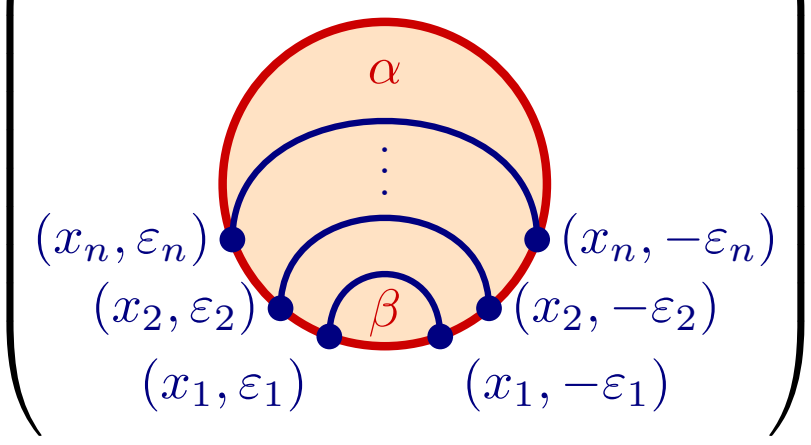

$\Bigg)(1) : {}^\dagger X \otimes X \longrightarrow 1_\alpha \quad (2.18)$

and

$\mathrm{coev}_X =$ 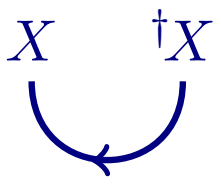

$= \mathcal{Z}\Bigg($ 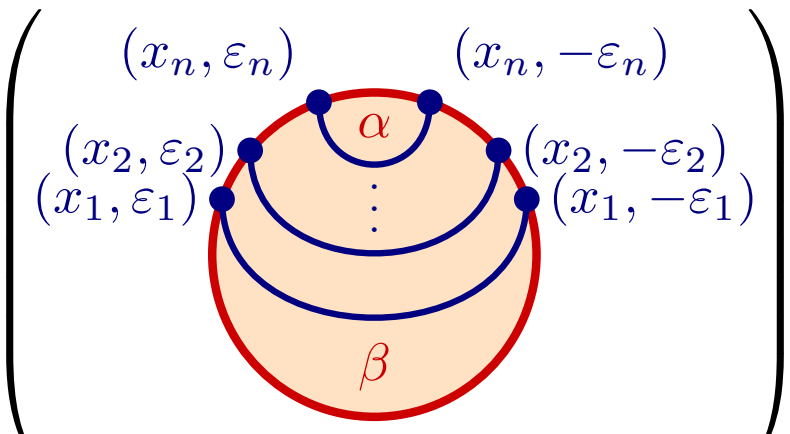

$\Bigg)(1) : 1_\beta \longrightarrow X \otimes {}^\dagger X\,. \quad (2.19)$

The maps $\widetilde{\mathrm{ev}}_X, \widetilde{\mathrm{coev}}_X$ exhibiting $X^\dagger$ as a right adjoint are defined analogously, by reversing all orientations and orders in (2.18) and (2.19). Proving that the Zorro moves hold is straighforward, for example

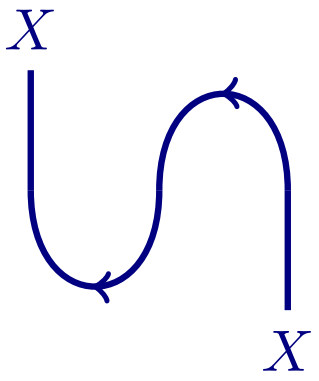

$= \mathcal{Z}\Bigg($ 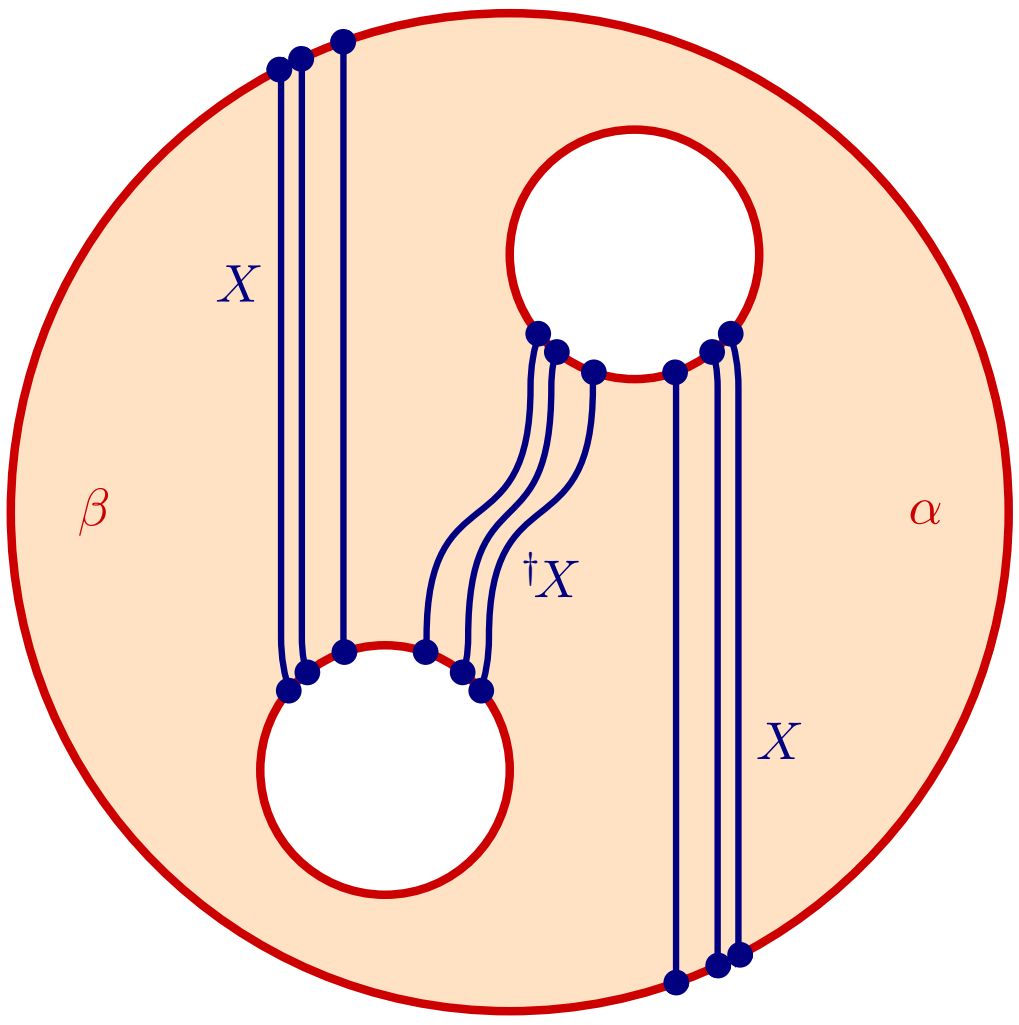

$\Bigg)(\mathrm{coev}_X \otimes_{\Bbbk} \mathrm{ev}_X)$

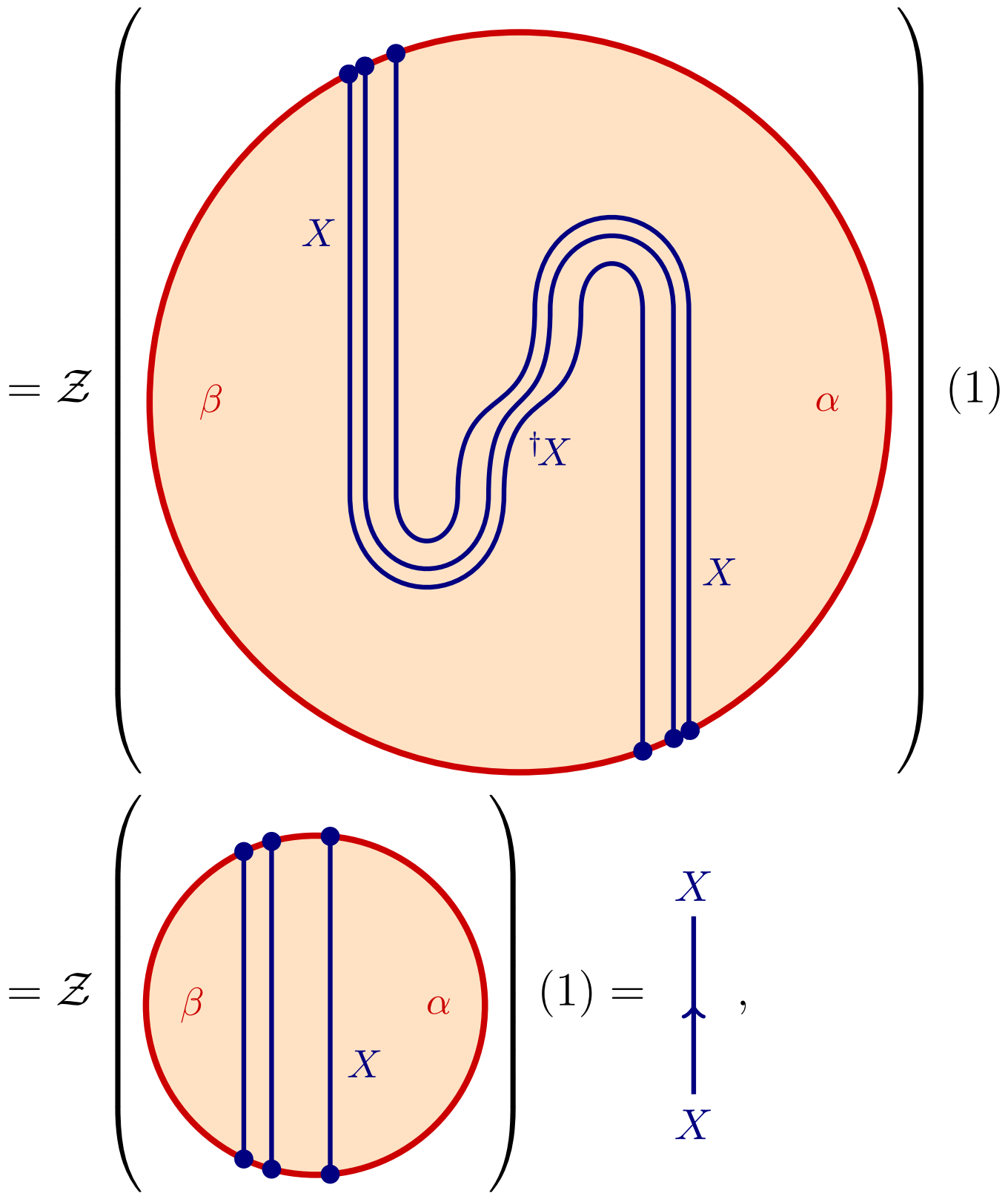

where in the second step we used functoriality of $\mathcal{Z}$, and in the third step we used isotopy invariance in $\mathrm{Bord}_2^{\mathrm{def}}(\mathbb{D})$. The pivotality conditions (2.7) are proved similarly, and we have arrived at the following result:

**Theorem 2.1.** Every 2-dimensional defect TQFT $\mathcal{Z}$ gives rise to a $\Bbbk$-linear pivotal 2-category $\mathcal{B}_{\mathcal{Z}}$ as constructed above.

By construction, string diagrams in $\mathcal{B}_{\mathcal{Z}}$ are correlators of $\mathcal{Z}$. More precisely, a priori the 2-morphisms in $\mathcal{B}_{\mathcal{Z}}$ only capture the action of $\mathcal{Z}$ on defect bordisms that can be embedded into the plane. But what about the other bordisms for which this in not possible, such as the sphere $S^2$? Below in Section 3.1 we will see how sphere correlators can be evaluated in $\mathcal{B}_{\mathcal{Z}}$ – under one additional assumption. Namely, in order to get from arbitrary defect bordisms $\Sigma$ to string diagrams, one can project $\Sigma$ onto a fixed plane, and generically this leads to a string diagram.[3] If this projection is not injective, then patches of different phases of $\Sigma$ will overlap on the plane. Hence one is led to a multiplicative structure on the set $D_2$ labelling phases, i. e. the objects of $\mathcal{B}_{\mathcal{Z}}$. More precisely, for the general procedure to consistently express defect correlators as string diagrams in $\mathcal{B}_{\mathcal{Z}}$, i. e. to produce the same result for every generic projection, $\mathcal{B}_{\mathcal{Z}}$ should be *monoidal*. This is in fact true of all the examples associated to defect TQFTs

[3] The situation here is similar to that of knots in $\mathbb{R}^3$ and their knot diagrams.

I know of, and it is natural to conjecture that 2-dimensional defect TQFTs are equivalent to $\Bbbk$-linear monoidal pivotal 2-categories.

It is still an open problem to classify defect TQFTs, by proving a theorem along the lines of the above conjecture. Thus the situation is different from that of closed TQFTs (which are equivalent to commutative Frobenius algebras), and that of open/closed TQFTs. The latter are basically equivalent to Calabi-Yau categories and the 'Cardy condition', the details of which we review in Section 3.2.

Before moving on to examples, we mention that the classification question is already settled for another, related type of 'enhanced' TQFTs. This however comes at the price of defining such TQFTs as higher functors between higher categories, contrary to the approach discussed above where a higher category is *constructed from* an ordinary functor. Indeed, a *2-1-0-extended TQFT* is a symmetric monoidal 2-functor $\mathrm{Bord}_{2,1,0} \to \mathrm{Alg}_\Bbbk$, where roughly $\mathrm{Bord}_{2,1,0}$ has points, lines, and 2-manifolds with corners (all oriented) as objects, 1-, and 2-morphisms, respectively, while $\mathrm{Alg}_\Bbbk$ consists of finite-dimensional $\Bbbk$-algebras, bimodules, and bimodule maps. It was shown in [SP] that such extended oriented TQFTs are equivalent to separable symmetric Frobenius algebras in $\mathrm{Alg}_\Bbbk$. This is precisely as predicted by the cobordism hypothesis, as homotopy fixed points of the trivial $\mathrm{SO}(2)$-action on fully dualisable objects in $\mathrm{Alg}_\Bbbk$ are separable symmetric Frobenius algebras [HSV].

## 2.4 Examples

In this section we sketch a number of pivotal 2-categories which are known (those in Sections 2.4.1, 2.4.7, and 2.4.8) or believed (those in Sections 2.4.2–2.4.6) to arise from defect TQFTs as described above. Presenting all details and the necessary background would inflate the length of the exposition exponentially. Hence we content ourselves with a rough account which mainly aims to subsume all examples under the common heading of defect TQFT, and to showcase the broad spectrum of interesting pivotal 2-categories.

In fact most of the examples are not strict 2-categories, but their weak cousins called *bicategories* for which horizontal composition is associative and unital only up to coherent isomorphisms [Ben]. Luckily, every (pivotal) bicategory is equivalent to a (pivotal) 2-category,[4] so we just as well may work with the bicategories that 'occur in nature'.

### 2.4.1 State sum models

Separable symmetric Frobenius $\Bbbk$-algebras also appear as special cases in defect TQFT. Indeed, they are the objects of a bicategory $\mathrm{ssFrob}_\Bbbk$, whose 1-morphisms from $A$ to $B$ are finite-dimensional $B$-$A$-bimodules, and 2-morphisms are bimodule maps. Horizontal composition of $M : A \to B$ and $N : B \to C$ is the tensor product $N \otimes_B M$ over the intermediate algebra, and the unit 1-morphism $1_A$ is $A$ viewed as a bimodule over itself.[5] The left and right adjoints of a 1-morphism $M$ are given by the dual bimodule $M^*$. Thanks to the natural isomorphism $M^{**} \cong M$, the $\Bbbk$-linear bicategory $\mathrm{ssFrob}_\Bbbk$ is also pivotal.

[4] An analogous maximal 'strictification' result does not hold for tricategories, which explains part of the richness of 3-dimensional TQFT.

It was shown in [DKR, Sect. 3] that $\mathrm{ssFrob}_\Bbbk$ is equivalent to the 2-category $\mathcal{B}_{\mathcal{Z}^{\mathrm{ss}}}$ associated to a special defect TQFT $\mathcal{Z}^{\mathrm{ss}}$. To wit,

$$\mathcal{Z}^{\mathrm{ss}} : \mathrm{Bord}_2^{\mathrm{def}}\big(D_1^{\mathrm{ss}}, D_2^{\mathrm{ss}}, s, t\big) \longrightarrow \mathrm{Vect}_\Bbbk$$

is a *state sum model*, generalising the state sum constructions of closed [FHK] and open/closed [LP2] 2-dimensional TQFTs. Here the defect data consist of the set $D_2^{\mathrm{ss}}$ of separable symmetric Frobenius algebras, the set $D_1^{\mathrm{ss}}$ of $B$-$A$-bimodules $M$ for all $A, B \in D_2^{\mathrm{ss}}$, and we have $s(M) = A$ and $t(M) = B$.

As in the closed and open/closed case, the construction of $\mathcal{Z}^{\mathrm{ss}}$ involves a choice of triangulation for objects and bordisms in $\mathrm{Bord}_2^{\mathrm{def}}(\mathbb{D}^{\mathrm{ss}})$, a decoration of every triangulation with algebraic data, and a projection procedure that ensures that the construction is independent of the choice of triangulation. In the case of $\mathcal{Z}^{\mathrm{ss}}$, additional technicalities arise from the compatibility of triangulation and defect lines. We refer to [DKR] for all details and only note that on objects, $\mathcal{Z}^{\mathrm{ss}}$ is defined by

$\mathcal{Z}^{\mathrm{ss}}\Bigg($ 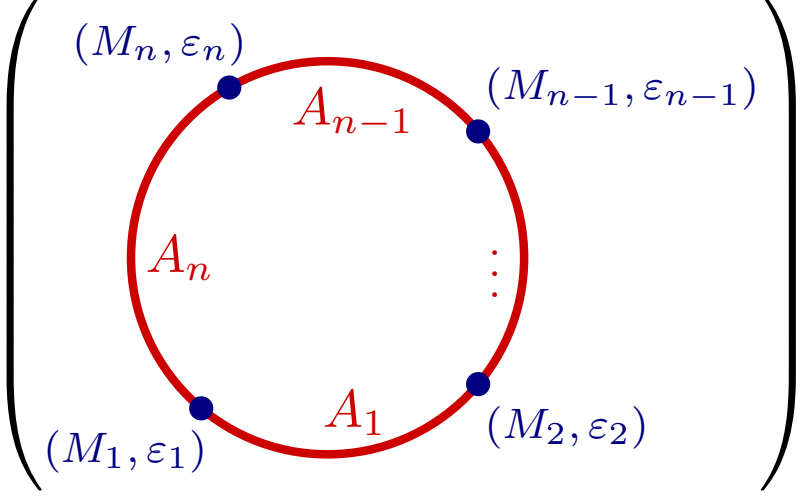

$\Bigg) = \circlearrowleft_{A_n}\Big(M_n^{\varepsilon_n} \otimes_{A_{n-1}} M_{n-1}^{\varepsilon_{n-1}} \otimes_{A_{n-2}} \cdots \otimes_{A_1} M_1^{\varepsilon_1}\Big)$

where $M^+ = M$ and $M^- = M^*$, and for an $A$-$A$-bimodule $M$, the vector space $\circlearrowleft_A M$ is defined to be the cokernel of the map $A \otimes M \to M$, $a \otimes m \mapsto am - ma$. It follows that $\mathcal{Z}^{\mathrm{ss}}(\textcircled{\scriptsize A}) = A/[A, A]$ is the 0-th Hochschild homology.

Finally we note that $\mathrm{ssFrob}_\Bbbk$ naturally has a monoidal structure, given by tensoring over the field $\Bbbk$. The trivial Frobenius algebra $\Bbbk$ is the unit object. Invertible objects are those algebras $A \in \mathrm{ssFrob}_\Bbbk$ for which there exist algebras $B, B' \in \mathrm{ssFrob}_\Bbbk$ such that $A \otimes_\Bbbk B$ and $B' \otimes_\Bbbk A$ are Morita equivalent to $\Bbbk$. It follows that isomorphism classes of invertible objects in $\mathrm{ssFrob}_\Bbbk$ precisely form the *Brauer group* of $\Bbbk$.

### 2.4.2 Algebraic geometry

Next we consider a more geometric example, the bicategory $\mathcal{V}ar$. It has smooth and proper varieties as objects, 1-morphisms are Fourier-Mukai kernels, and 2-

[5] It follows that invertible 1-morphisms in $\mathrm{ssFrob}_\Bbbk$ are precisely Morita equivalences.

morphisms are their maps up to quasi-isomorphism. Hence for $U, V \in \mathcal{V}ar$, we have that

$$\mathcal{V}ar(U, V) = \mathbb{D}^{\mathrm{b}}(\mathrm{coh}(U \times V)) =: \mathbb{D}(U \times V)$$

is the bounded derived category of coherent sheaves on the product space. Horizontal composition of kernels $\mathcal{E}, \mathcal{F}$ is the convolution $\mathcal{E} \circ \mathcal{F}$, and the unit $1_U$ is the structure sheaf $\mathcal{O}_{\Delta_U}$ of the diagonal $\Delta_U \subset U \times U$.

Adjunctions in $\mathcal{V}ar$ were studied in detail in [CW]. The 'naive' adjoint of $\mathcal{E} \in \mathcal{V}ar(U, V)$ is $\mathcal{E}^\vee := \underline{\mathrm{Hom}}_{\mathbb{D}(U \times V)}(\mathcal{E}, U \times V)$, but by Grothendieck duality the true adjoints are obtained by 'twisting with the Serre kernel': we have

$$^\dagger\mathcal{E} = \mathcal{E}^\vee \circ \Sigma_V \, , \quad \mathcal{E}^\dagger = \Sigma_U \circ \mathcal{E}^\vee \, ,$$

where $\Sigma_U$ is obtained from the canonical line bundle $\omega_U$ as $\Sigma_U = (\Delta_U)_* \omega_U[\dim U]$. It follows that $\mathcal{V}ar$ is not pivotal 'on the nose', but only has *polite dualities* as explained in [CW].

One way to look at the failure of $\mathcal{V}ar$ being strictly pivotal is via its interpretation in terms of *B-twisted sigma models*, whose defects were also studied in [Sar, AS]. Indeed, in theoretical physics to every object $U$ in $\mathcal{V}ar$ one associates a field theory called an '$\mathcal{N} = (2, 2)$ supersymmetric sigma model'. There is a procedure called 'topological B-twist' [HKK+, Ch. 16] that produces a closed 2-dimensional TQFT from $U$ if it is a Calabi-Yau variety, i. e. $\omega_U$ is trivial. Indeed, in this case Hochschild homology and cohomology coincide up to shift and are isomorphic to Dolbeault cohomology $H_{\bar\partial}(U)$. Then together with the Mukai pairing $H_{\bar\partial}(U)$ is a commutative Frobenius algebra in the category of graded vector spaces $\mathrm{Vect}^{\mathrm{gr}}_{\mathbb{C}}$. Hence $H_{\bar\partial}(U)$ describes a closed TQFT $\mathrm{Bord}_2 \to \mathrm{Vect}^{\mathrm{gr}}_{\mathbb{C}}$. The fact that $H_{\bar\partial}(U)$ is only graded commutative can be traced back to the supersymmetry of the original sigma model.

To $\mathcal{V}ar$ one can associate the defect data $\mathbb{D}^{\mathrm{B}}$ given by $\mathbb{D}^{\mathrm{B}}_2 = \mathrm{Obj}(\mathcal{V}ar)$, $\mathbb{D}^{\mathrm{B}}_1 = \{\mathrm{Obj}(\mathcal{V}ar(U, V))\}_{U,V \in \mathbb{D}^{\mathrm{B}}_2}$ and obvious source and target maps. Then it is natural to conjecture that there is a defect TQFT $\mathcal{Z}^{\mathrm{B}} : \mathrm{Bord}^{\mathrm{def}}_2(\mathbb{D}^{\mathrm{B}}) \to \mathrm{Vect}^{\mathrm{gr}}_{\mathbb{C}}$ whose associated 2-category $\mathcal{B}_{\mathcal{Z}^{\mathrm{B}}}$ is equivalent to $\mathcal{V}ar$. Furthermore, taking products of varieties gives $\mathcal{V}ar$ a monoidal structure.

### 2.4.3 Symplectic geometry

Another geometric example is the pivotal bicategory Symp, which is "dual" to $\mathcal{V}ar$ in the sense that it is believed to arise from a defect TQFT that collects all *A-twisted sigma models*. (For the subbicategories of Calabi-Yau varieties this relation is expected to be a generalisation of mirror symmetry.) The bicategory Symp is discussed in detail in the review article [Weh, Sect. 3.5]. Objects are symplectic manifolds $M \equiv (M, \omega)$, and one writes $M^- = (M, -\omega)$ for the reversed symplectic structure. 1-morphisms $M \to N$ are Lagrangian correspondences $\underline{L}$,

i. e. chains of Lagrangian submanifolds $L_{i,j} \subset M_i^- \times M_j$ of the form

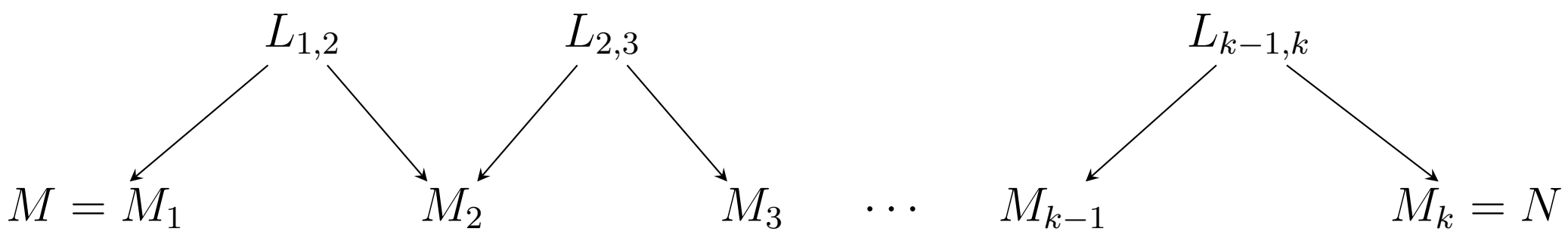

while 2-morphisms form the quilted Floer homology groups $\mathrm{Hom}(\underline{L}, \underline{L}') = HF(\underline{L}, \underline{L}')$. Hence the categories of 1-morphisms in Symp are Donaldson-Fukaya categories.

The left and right adjoint of $\underline{L} = (L_{1,2}, L_{2,3}, \dots, L_{k-1,k})$ is the reversed correspondence $\underline{L}^T = (L^T_{k-1,k}, L^T_{k-2,k-1}, \dots, L^T_{1,2})$, where $L^T_{i,j} \subset M_j^- \times M_i$ is the image of $L_{i,j}$ under the transposition $M_i \times M_j \to M_j \times M_i$. Taking products of symplectic manifolds makes Symp a monoidal 2-category.

### 2.4.4 Landau-Ginzburg models

That left and right adjoints agree only up to a certain twist is a common phenomenon, which was formalised and called 'graded pivotal' in [CM2, Def. 7.1]. An example of such a graded pivotal bicategory is that of (affine) *Landau-Ginzburg models* $\mathcal{LG}$. Its objects are isolated singularities, i. e. polynomials $W \in \mathbb{C}[x_1, \dots, x_n]$ for some $n \in \mathbb{N}$ such that the Jacobi ring

$$\mathrm{Jac}_W := \mathbb{C}[x_1, \dots, x_n]/(\partial_{x_1} W, \dots, \partial_{x_n} W)$$

is finite-dimensional over $\mathbb{C}$.[6] A 1-morphism from $W \in \mathbb{C}[x_1, \dots, x_n] \equiv \mathbb{C}[x]$ to $V \in \mathbb{C}[z_1, \dots, z_m] \equiv \mathbb{C}[z]$ is a *matrix factorisation* $X$ of $V - W$. This means that $X$ is a free finite-rank $\mathbb{Z}_2$-graded $\mathbb{C}[z,x]$-module $X = X^0 \oplus X^1$ together with an odd map $d_X \in \mathrm{End}^1_{\mathbb{C}[z,x]}(X)$ such that $d_X^2 = (V - W) \cdot 1_X$. A 2-morphism between $X, Y \in \mathcal{LG}(W,V)$ is an even $\mathbb{C}[z,x]$-linear map $\Phi : X \to Y$ up to homotopy with respect to the twisted differentials $d_X$ and $d_Y$. If one chooses bases of $X$ and $Y$, then $d_X, d_Y$ and $\Phi$ are represented by odd and even matrices, respectively.

Horizontal composition in $\mathcal{LG}$ is given by tensoring over the intermediate polynomial ring.[7] The unit $1_W$ is a deformation of the Koszul complex of $(\partial_{x_1} W, \dots, \partial_{x_n} W)$, which in the simplest case of $W = x^d$ means that $d_{1_W}$ is represented by the matrix $\left(\begin{smallmatrix} 0 & x-x' \\ (x^d - x'^d)/(x-x') & 0 \end{smallmatrix}\right)$. In general one finds that $\mathrm{End}(1_W) \cong \mathrm{Jac}_W$.

Not only horizontal composition, but also adjunctions in $\mathcal{LG}$ are under very good control. Up to a shift, ${}^\dagger X$ and $X^\dagger$ are simply given by $\mathrm{Hom}_{\mathbb{C}[z,x]}(X, \mathbb{C}[z,x])$, but also the adjunction maps $\mathrm{ev}_X, \mathrm{coev}_X$ etc. were computed explicitly in terms of Atiyah classes in [CM2]. In this way we could show that $\mathcal{LG}$ is graded pivotal. Furthermore, one obtains neat formulas for quantum dimension and traces (recall (2.8) and (2.9)), for example

$$\dim_{\mathrm{r}}(X) = (-1)^{\binom{m+1}{2}} \operatorname{Res}_{\mathbb{C}[x,z]/\mathbb{C}[z]} \left[ \frac{\operatorname{str}\left(\partial_{x_1} d_X \dots \partial_{x_n} d_X \, \partial_{z_1} d_X \dots \partial_{z_m} d_X\right) \underline{\mathrm{d}x}}{\partial_{x_1} W \dots \partial_{x_n} W} \right]$$

where $X \in \mathcal{LG}(W, V)$ is as above.

[6] In fact there is a graded pivotal bicategory $\mathcal{LG}_\Bbbk$ for any commutative ring $\Bbbk$ as explained in [CM2], but the finiteness condition on its objects is more involved for $\Bbbk \neq \mathbb{C}$.

[7] This can be algorithmically computed by splitting an idempotent [DM], which was used in [CM1] to compute the homological knot invariants of Khovanov and Rozansky [KR].

Using the construction of [Shu] one may verify that $\mathcal{LG}$ is a monoidal bicategory, where the tensor product of $W \in \mathbb{C}[x]$ with $V \in \mathbb{C}[z]$ is $W + V \in \mathbb{C}[z, x]$, and $0 \in \mathbb{C}$ is the monoidal unit. It is expected (but not proven) that there exists a defect TQFT $\mathcal{Z}^{\mathrm{LG}} : \mathrm{Bord}_2^{\mathrm{def}}(\mathbb{D}^{\mathrm{LG}}) \to \mathrm{Vect}_{\Bbbk}^{\mathbb{Z}_2}$ such that $\mathcal{LG}_{\Bbbk}$ is equivalent to the 2-category $\mathcal{B}_{\mathcal{Z}^{\mathrm{LG}}}$ associated to $\mathcal{Z}^{\mathrm{LG}}$.

### 2.4.5 Differential graded categories

Twisted sigma models and Landau-Ginzburg models fit into a larger framework of differential graded (dg) categories [Toë]. Indeed, there is a bicategory $\mathcal{DG}_{\Bbbk}^{\mathrm{sat}}$ whose objects are 'saturated' dg categories, and whose 1- and 2-morphisms are certain resolutions of dg functors and natural transformations, respectively. Every object $U \in \mathcal{V}ar$, $M \in \mathrm{Symp}$, or $W \in \mathcal{LG}$ can be viewed as an object in $\mathcal{DG}_{\Bbbk}^{\mathrm{sat}}$ by taking the unique dg enhancemens of $\mathbb{D}(U)$, $\mathrm{Symp}(\mathrm{pt}, M)$, or $\mathcal{LG}(0, W)$, respectively. Then the bicategories $\mathcal{V}ar$, Symp and $\mathcal{LG}$ are quasi-equivalent to the (distinct) corresponding full subbicategories of $\mathcal{DG}_{\Bbbk}^{\mathrm{sat}}$. More generally, one might think of $\mathcal{DG}_{\Bbbk}^{\mathrm{sat}}$ as the bicategory of TQFTs arising from topologically twisting $\mathcal{N} = (2, 2)$ supersymmetric quantum field theories.

Working in the enlarged framework of $\mathcal{DG}_{\Bbbk}^{\mathrm{sat}}$ also has the advantage of comparing sigma models and Landau-Ginzburg models in a more natural context. In particular, homological mirror symmetry is a statement internal to $\mathcal{DG}_{\Bbbk}^{\mathrm{sat}}$, and it is tempting to speculate that there is a truly 2-categorical generalisation of mirror symmetry, taking place in a suitable enlargement of $\mathcal{DG}_{\Bbbk}^{\mathrm{sat}}$.[8]

As shown in [BFK, App. A.2], $\mathcal{DG}_{\Bbbk}^{\mathrm{sat}}$ is equivalent to the more manageable bicategory $\mathcal{DG}_{\Bbbk}^{\mathrm{sp}}$ of smooth and proper dg algebras. The 1-morphisms $M : A \to B$ in $\mathcal{DG}_{\Bbbk}^{\mathrm{sp}}$ are perfect dg $(A^{\mathrm{op}} \otimes_{\Bbbk} B)$-bimodules, i. e. $\mathrm{Hom}_{\mathbb{D}(A^{\mathrm{op}} \otimes_{\Bbbk} B)}(M, -)$ commutes with arbitrary coproducts, and 2-morphisms are maps of dg bimodules up to quasi-isomorphisms. Horizontal composition is the left-derived tensor product over the intermediate algebra. It follows from the discussion in [BFK] that $\mathcal{DG}_{\Bbbk}^{\mathrm{sp}}$ is graded pivotal, and analogously to the situation in Sections 2.4.2 and 2.4.4, adjoints are given by the naive dual together with a twist by Serre functors.

Tensoring over the field $\Bbbk$ gives the bicategories $\mathcal{DG}_{\Bbbk}^{\mathrm{sat}}$ and $\mathcal{DG}_{\Bbbk}^{\mathrm{sp}}$ a natural monoidal structure [Toë].

[8] The process of orbifold completion discussed in Section 2.4.8 below may play a role here.

### 2.4.6 Categorified quantum groups

Pivotal 2-categories also feature prominently in higher representation theory (see e. g. [Lau2, Sect. 1] for a gentle introduction), which in turn plays a unifying role in the theory of homological link invariants. A key idea behind '2-Kac-Moody algebras' [Rou] and 'categorified quantum groups' [Lau1, KL] is to represent them not on vector spaces, but on linear categories. This leads one to replace the (quantum) Serre relations, which are equalities between expressions involving the generators $E_i, F_j$, by natural transformations $\eta_k$ between the corresponding expressions of functors $\mathcal{E}_i, \mathcal{F}_j$. This theory is considerably richer than its classical counterpart, partly because of nontrivial relations which have to be imposed on the $\eta_k$.

Associated to any Kac-Moody algebra $\mathfrak{g}$ there is a $\Bbbk$-linear pivotal 2-category $\mathcal{U}_Q(\mathfrak{g})$. The precise definition fills several pages (cf. the above references), but here is a sketch: First one picks a 'Cartan datum and choice of scalars $Q$'; this in particular gives a weight lattice $X$ with simple roots $\alpha_i$, and a symmetrisable generalised Cartan matrix. Objects of $\mathcal{U}_Q(\mathfrak{g})$ are simply weights $\lambda \in X$, and 1-morphisms are formal polynomials in expressions of the form

$$\mathbf{1}_\lambda \, , \quad \mathbf{1}_{\lambda+\alpha_i}\mathcal{E}_i = \mathbf{1}_{\lambda+\alpha_i}\mathcal{E}_i\mathbf{1}_\lambda = \mathcal{E}_i\mathbf{1}_\lambda \, , \quad \mathbf{1}_{\lambda-\alpha_i}\mathcal{F}_i = \mathbf{1}_{\lambda-\alpha_i}\mathcal{F}_i\mathbf{1}_\lambda = \mathcal{F}_i\mathbf{1}_\lambda \, .$$

2-morphisms are $\Bbbk$-spans of compositions of certain string diagrams which encode the categorification of the Serre relations for $\mathfrak{g}$, subject to a list of relations. These relations in particular describe biadjunctions (again up to shifts) between $\mathcal{E}_i\mathbf{1}_\lambda$ and $\mathcal{F}_i\mathbf{1}_\lambda$. As explained in [BHLW], the parameters $Q$ can be chosen such that these adjunctions become a strictly pivotal structure on $\mathcal{U}_Q(\mathfrak{g})$. It is an open problem to determine whether there is a natural monoidal structure on $\mathcal{U}_Q(\mathfrak{g})$.

### 2.4.7 Surface defects in 3-dimensional TQFT

One generally expects $n$-dimensional TQFTs to appear as defects of codimension 1 in $(n+1)$-dimensional TQFTs. At least for $n = 2$ this is a rigorous result also from the algebraic perspective: In [CMS] we introduced 3-dimensional defect TQFTs $\mathcal{Z}$, from which we went on to construct a certain type of 3-category $\mathcal{T}_{\mathcal{Z}}$. More precisely, $\mathcal{T}_{\mathcal{Z}}$ is a $\Bbbk$-linear 'Gray category with duals' – categorifying the construction of Section 2.3. This implies in particular that $\mathcal{T}_{\mathcal{Z}}(u, v)$ is a $\Bbbk$-linear pivotal 2-category for all $u, v \in \mathcal{T}_{\mathcal{Z}}$. And since the objects (interpreted as 'surface defects') of $\mathcal{T}_{\mathcal{Z}}(u, v)$ are the 1-morphisms of a 3-category, $\mathcal{T}_{\mathcal{Z}}(u, v)$ is monoidal whenever $u = v$.

### 2.4.8 Orbifold completion

As a final source of pivotal 2-categories we point to the procedure of 'orbifold completion' of [CR]. Inspired by the orbifold construction from finite group actions and their generalisation in rational conformal field theory [FFRS], orbifold completion takes a pivotal bicategory $\mathcal{B}$ as input and produces a new pivotal bicategory $\mathcal{B}_{\text{orb}}$, into which $\mathcal{B}$ fully embeds. It is a completion because there is an equivalence $(\mathcal{B}_{\text{orb}})_{\text{orb}} \cong \mathcal{B}_{\text{orb}}$.

Objects of $\mathcal{B}_{\text{orb}}$ are pairs $(\alpha, A)$ where $\alpha \in \mathcal{B}$ and $A$ is a separable symmetric Frobenius algebra – not necessarily in $\text{Vect}_{\Bbbk}$, but in the category $\mathcal{B}(\alpha, \alpha)$. A 1-morphism $(\alpha, A) \to (\beta, B)$ in $\mathcal{B}_{\text{orb}}$ is a 1-morphism $X \in \mathcal{B}(\alpha, \beta)$ together with the structure of a $B$-$A$-bimodule, and 2-morphisms are those in $\mathcal{B}$ which are also bimodule maps. Horizontal composition is the tensor product over the intermediate algebra, and $1_{(\alpha,A)} = A$ viewed as an $A$-$A$-bimodule.

The simplest case of orbifold completion reproduces state sum models. Here as the input bicategory $\mathcal{B}$ one takes the 'trivial' bicategory $B\text{Vect}_{\Bbbk}$ which has a single object with $\text{Vect}_{\Bbbk}$ as its endomorphism category. Then by construction and in the notation of Section 2.4.1 we have

$$\text{ssFrob}_{\Bbbk} \cong (B\text{Vect}_{\Bbbk})_{\text{orb}} \,.$$

Examples of separable Frobenius algebras which do not live in $\text{Vect}_{\Bbbk}$ are provided by group actions. For a finite group $G$ and some pivotal bicategory $B$, let $D_g \in \mathcal{B}(\alpha, \alpha)$ for all $g \in G$ such that $D_g \otimes D_h \cong D_{gh}$ coherently, and let $\mathcal{B}(\alpha, \alpha)$ have finite sums. Then there are as many inequivalent separable Frobenius algebras structures on $A_G := \bigoplus_{g \in G} D_g$ as there are elements in $H^2(G, \Bbbk^\times)$ [BCP2], and $A_G$ is symmetric if its Nakayama automorphism is the identity. Interestingly, not all separable Frobenius algebras come from group actions, cf. [CRCR].

The adjoint of $X \in \mathcal{B}_{\text{orb}}((\alpha, A), (\beta, B))$ is ${}^\dagger X = X^\dagger \in \mathcal{B}(\beta, \alpha)$ together with the adjunction maps

$\text{ev}_X =$ 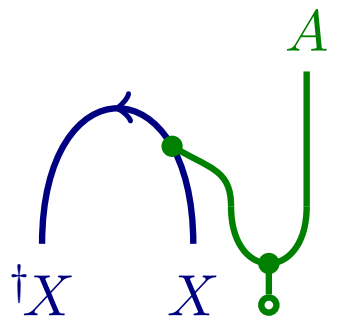 $\circ\, \xi \,, \quad \text{coev}_X = \vartheta \,\circ$ 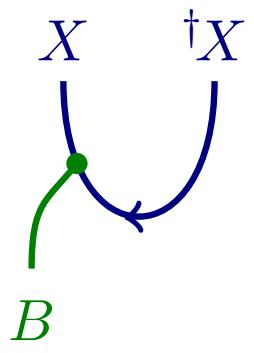

where $\xi : {}^\dagger X \otimes_B X \to {}^\dagger X \otimes X$ and $\vartheta : X \otimes {}^\dagger X \to X \otimes_A {}^\dagger X$ are the splitting and projection maps (which we require to exist, cf. [CR, Lem. 2.3]). In fact, this continues to hold if $A$ and $B$ are not required to be symmetric, but then their actions on ${}^\dagger X$ and $X^\dagger$ are twisted by Nakayama automorphisms as explained in [CR, Sec. 4.3], generalising the situation with Serre functors in Sections 2.4.2, 2.4.4 and 2.4.5, see [BCP1, CQV].

If $\mathcal{B}$ is monoidal then $\mathcal{B}_{\text{orb}}$ is expected to be monoidal as well. One may either argue along the lines of [Shu], or via a universal property for the operation $(-)_{\text{orb}}$.

# 3 Open/closed TQFTs from pivotal 2-categories

In Section 2.1.1 we easily obtained closed and open/closed TQFTs from defect TQFTs, simply by forgetting part of the structure of the defect bordism category $\mathrm{Bord}_2^{\mathrm{def}}(\mathbb{D})$. In the present section we elucidate this relation by constructing an open/closed TQFT from any object in a $\Bbbk$-linear pivotal bicategory (satisfying two natural assumptions). We treat the purely closed case in Section 3.1, and the fully open/closed case in Section 3.3. The algebraic description of open/closed TQFTs in terms of Calabi-Yau categories is reviewed in Section 3.2.

## 3.1 Closed TQFTs from pivotal 2-categories

For completeness, we briefly recall how closed TQFTs $\mathcal{Z}^{\mathrm{c}} : \mathrm{Bord}_2 \to \mathrm{Vect}_{\Bbbk}$ are equivalent to commutative Frobenius algebras in $\mathrm{Vect}_{\Bbbk}$. By a classical result [Koc], the bordism category $\mathrm{Bord}_2$ is generated by

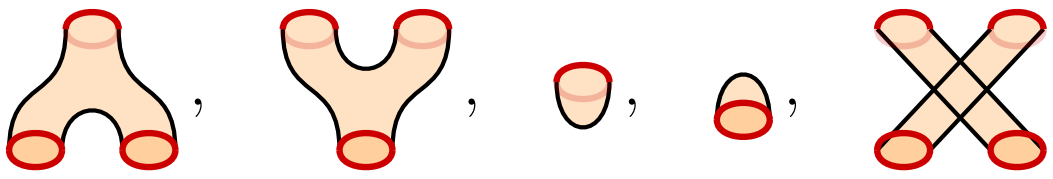 (3.1)

where from now on in-coming/out-going boundaries are placed at the bottom/top of pictures of bordisms. The relations between the generators (3.1) precisely say that the vector space $\mathcal{Z}^{\mathrm{c}}(S^1)$ together with multiplication $\mathcal{Z}^{\mathrm{c}}(\;)$, unit $\mathcal{Z}^{\mathrm{c}}(\;)(1)$ and pairing $\mathcal{Z}^{\mathrm{c}}(\;) \circ \mathcal{Z}^{\mathrm{c}}(\;)$ is a commutative Frobenius algebra.

Given a $\Bbbk$-linear pivotal bicategory $\mathcal{B}$, we would now like to construct a commutative Frobenius algebra $A_\alpha$ naturally associated to every object $\alpha \in \mathcal{B}$. Our interpretation (2.10) of such bicategories in the context of defect TQFT suggests that the underlying vector space of $A_\alpha$ should be the space of endomorphisms of the unit $1_\alpha$:

$$A_\alpha = \mathrm{End}(1_\alpha)\,.$$

Indeed, we interpret $1_\alpha$ as the invisible defect, and elements $\phi \in \mathrm{End}(1_\alpha)$ correspond to operators living on an invisible line with no further 'defect conditions':

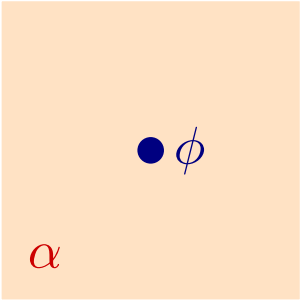

$\in \mathrm{End}(1_\alpha)\,.$

Hence in standard jargon, $\phi$ is a field 'inserted in the bulk' of the 'theory' $\alpha$.

As the 2-endomorphism space in a $\Bbbk$-linear bicategory, $A_\alpha$ is manifestly an associative unital $\Bbbk$-algebra. Furthermore, using the interchange law (2.3) one can show that in any monoidal category the endomorphisms of the unit form a commutative monoid [EGNO, Prop. 2.2.10].

It remains to endow $A_\alpha$ with a nondegenerate pairing which is compatible with multiplication. There are various ways to ensure the existence of such a pairing. We will see how it follows from

**Assumption 3.1.** The $\Bbbk$-linear pivotal bicategory $\mathcal{B}$ is monoidal[9] with duals such that $\mathcal{B}(\mathbb{O}, \mathbb{O}) \cong \mathrm{Vect}_\Bbbk$, where $\mathbb{O} \in \mathcal{B}$ is the unit object.

We observe that all our examples in Section 2.4 are known or expected to satisfy this condition. Further we recall that as noted after Theorem 2.1 it is natural to expect the bicategories arising from defect TQFTs to come with a monoidal structure: for $\alpha, \beta \in \mathcal{B}$ their tensor product $\alpha \,\square\, \beta$ corresponds to the 'tensor product theory' attached to the fusion of two bordism patches labelled $\alpha$ and $\beta$. The unit $\mathbb{O} \in \mathcal{B}$ corresponds to the 'trivial theory'. This is consistent with $A_\mathbb{O} = \mathrm{End}(1_\mathbb{O}) \cong \Bbbk$ being the trivial commutative Frobenius algebras in $\mathrm{Vect}_\Bbbk \cong \mathcal{B}(\mathbb{O}, \mathbb{O})$. Finally, the dual $\alpha^\#$ of an object $\alpha \in \mathcal{B}$ is interpreted as a label for a bordism patch with opposite orientation.

How can we endow $A_\alpha = \mathrm{End}(1_\alpha)$ with a nondegenerate pairing using Assumption 3.1? Recall that the pairing of the Frobenius algebra associated to a closed TQFT $\mathcal{Z}^{\mathrm{c}}$ is the 'sphere correlator'

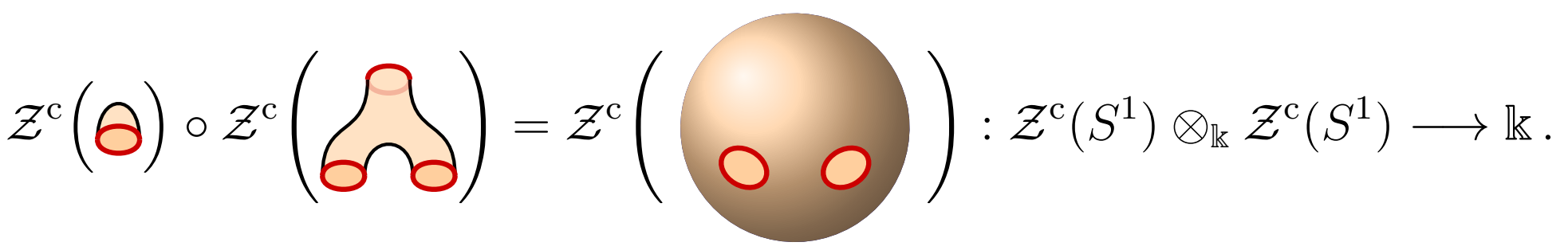

From a bicategory $\mathcal{B}$ as in Assumption 3.1 we can build such pairings by mimicking the above construction as follows. We consider a sphere $S^2$ labelled by $\alpha \in \mathcal{B}$. Then we project the sphere onto some plane, producing a disc:

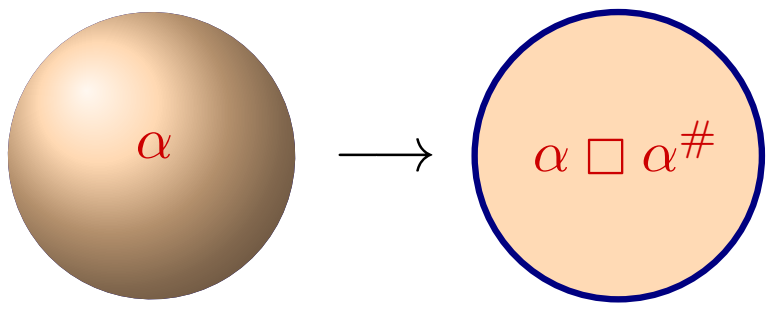

The disc is labelled by the product $\alpha \,\square\, \alpha^\#$ because the rear part of the sphere has opposite orientation with respect to the plane.

We now interpret the disc as a string diagram in $\mathcal{B}$! The preimage of the disc's boundary is simply the great circle on the sphere which is parallel to the chosen projection plane. Since any plane will do, this great circle is nothing special – in fact it is invisible on the sphere. Hence we label it with the 'invisible defect' $1_\alpha \in \mathcal{B}(\alpha, \alpha)$.

Labelling the inside and outside of the sphere with the trivial object $\mathbb{O}$, the boundary of the disc is labelled by the 1-morphism

$$\widetilde{1}_\alpha \in \mathcal{B}\big(\mathbb{O}, \alpha \,\square\, \alpha^\#\big)$$

[9] Monoidal bicategories are defined for example in [SP] or [Shu]. By a 'monoidal bicategory with duals' we mean a 'Gray category with duals and only a single object' as defined e. g. in [BMS, Sect. 3.3] or [CMS, Sect. 3.2.2].

corresponding to $1_\alpha$ under duality in $\mathcal{B}$. So finally the 'sphere correlator' for $A_\alpha = \mathrm{End}(1_\alpha)$ is defined to be

$$\big\langle -,- \big\rangle_\alpha : A_\alpha \otimes_{\Bbbk} A_\alpha \longrightarrow \Bbbk\,, \quad \phi_1 \otimes \phi_2 \longmapsto$$ 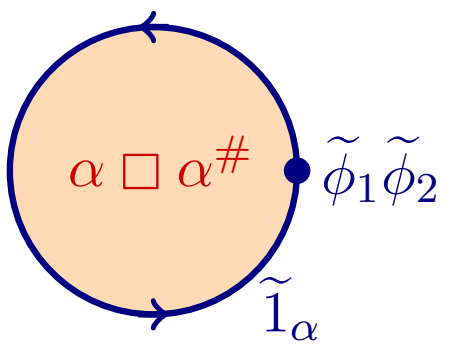

$$\in \mathrm{End}(1_{\mathbb{O}}) \cong \Bbbk \quad (3.2)$$

where the string diagram represents the map $\mathrm{ev}_{\widetilde{1}_\alpha} \circ (1_{\widetilde{1}_\alpha} \otimes \widetilde{\phi}_1 \widetilde{\phi}_2) \circ \widetilde{\mathrm{coev}}_{\widetilde{1}_\alpha} : \Bbbk \to \Bbbk$, i. e. the trace $\mathrm{tr}_1(\widetilde{\phi}_1 \widetilde{\phi}_2)$. Thus the pairing $\langle -,- \rangle_\alpha$ is manifestly compatible with multiplication, and it follows from Theorem 3.6 below that $\langle -,- \rangle_\alpha$ is nondegenerate. Hence we have arrived at a closed TQFT:

**Theorem 3.2.** Let $\mathcal{B}$ be a bicategory satisfying Assumption 3.1. Then for every $\alpha \in \mathcal{B}$, the vector space $\mathrm{End}(1_\alpha)$ naturally has the structure of a commutative Frobenius algebra.

## 3.2 Open/closed TQFTs and Calabi-Yau categories

An *open/closed* TQFT [Laz, MS, LP1] is a symmetric monoidal functor

$$\mathcal{Z}^{\mathrm{oc}} : \mathrm{Bord}_2^{\mathrm{oc}}(B) \longrightarrow \mathrm{Vect}_{\Bbbk} \quad (3.3)$$

where $B$ is some set, whose elements are referred to as *boundary conditions*. The objects of the bordism category $\mathrm{Bord}_2^{\mathrm{oc}}(B)$ are disjoint unions of circles $S^1$ and unit intervals $I_{ab}$ whose endpoints are labelled by $a, b \in B$:

$I_{ab} =$ 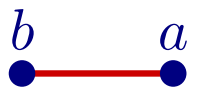 .

Morphisms in $\mathrm{Bord}_2^{\mathrm{oc}}(B)$ are bordism classes generated by the list (3.1) together with the classes of the following decorated manifolds with corners for all $a, b, c, d \in B$:

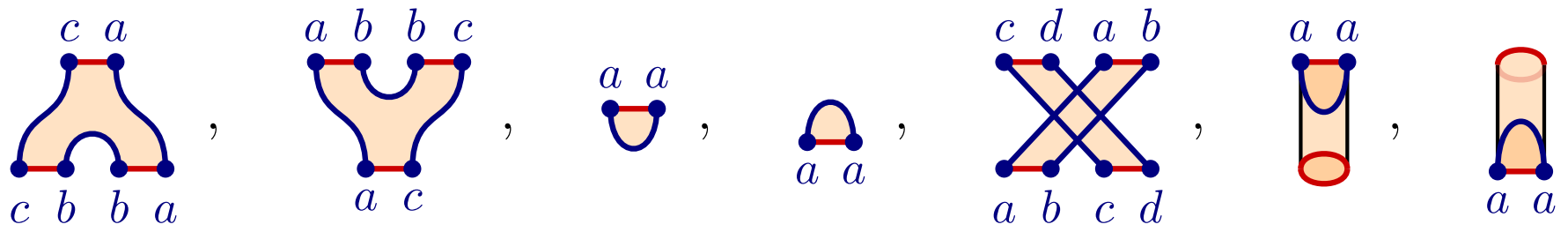

(3.4)

subject to the relations (3.5)–(3.12) below as well as the relations for 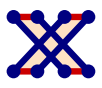 which state that together with the twist 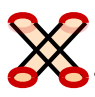, the category $\mathrm{Bord}_2^{\mathrm{oc}}(B)$ has a symmetric monoidal structure. It follows that $\mathrm{Bord}_2$ is a non-full subcategory of $\mathrm{Bord}_2^{\mathrm{oc}}(B)$.

The relations on the generators (3.4) are the intuitively clear equalities

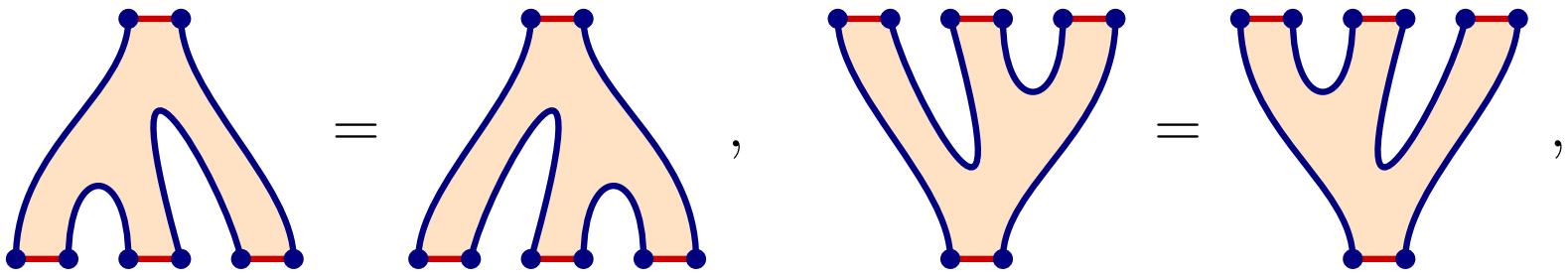

(3.5)

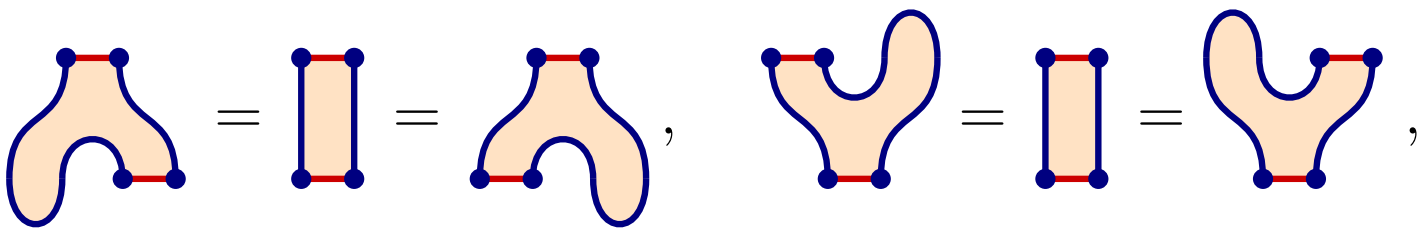

(3.6)

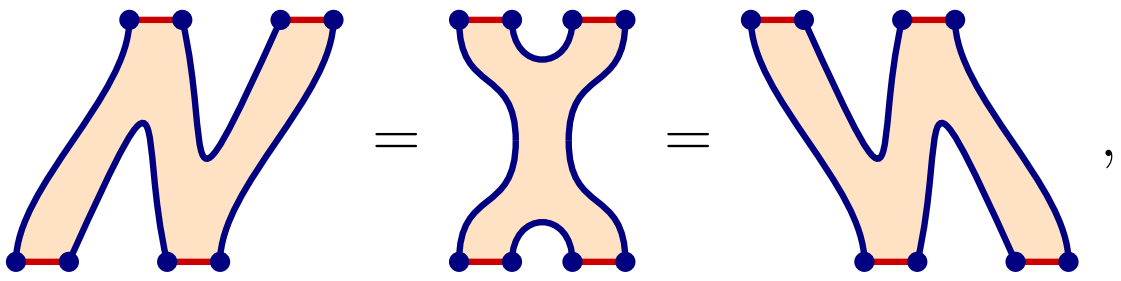

(3.7)

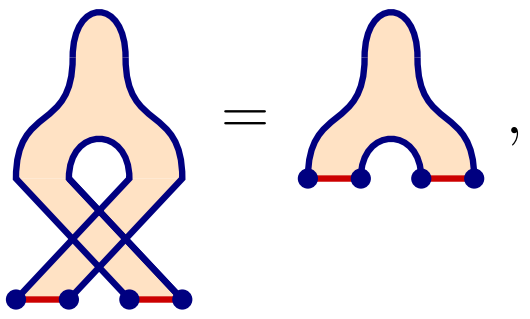

(3.8)

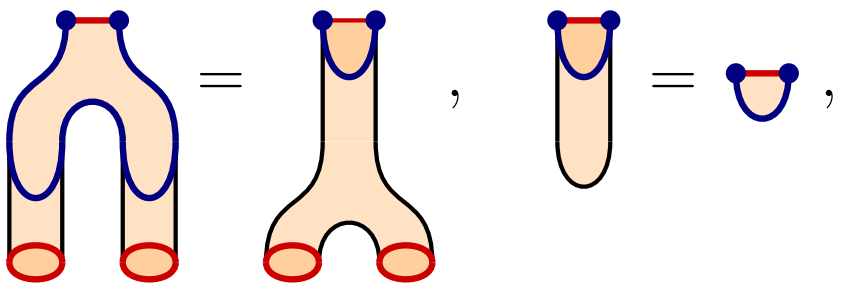

(3.9)

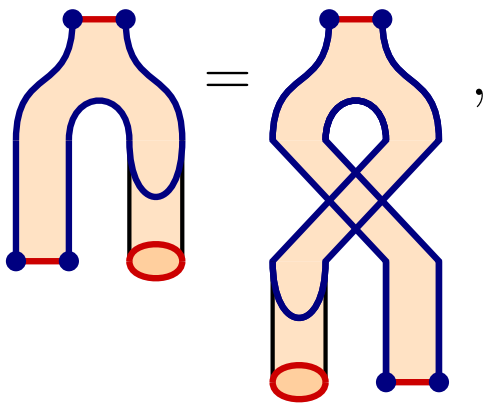

(3.10)

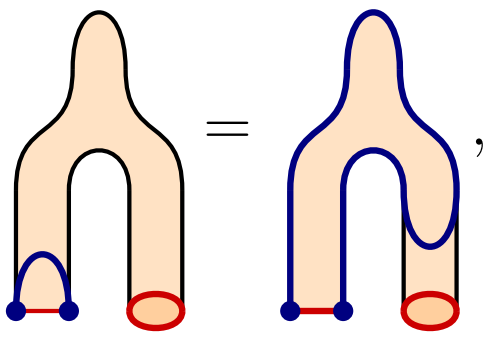

(3.11)

as well as

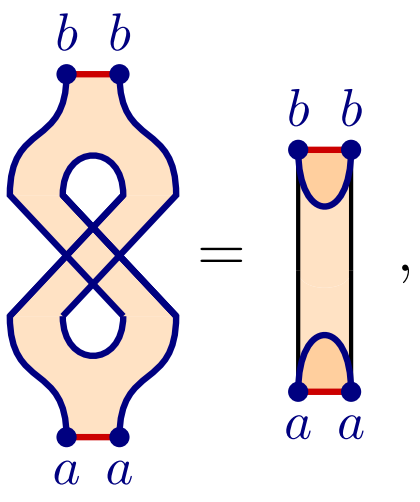

(3.12)

which can be understood as the sequence of diffeomorphisms

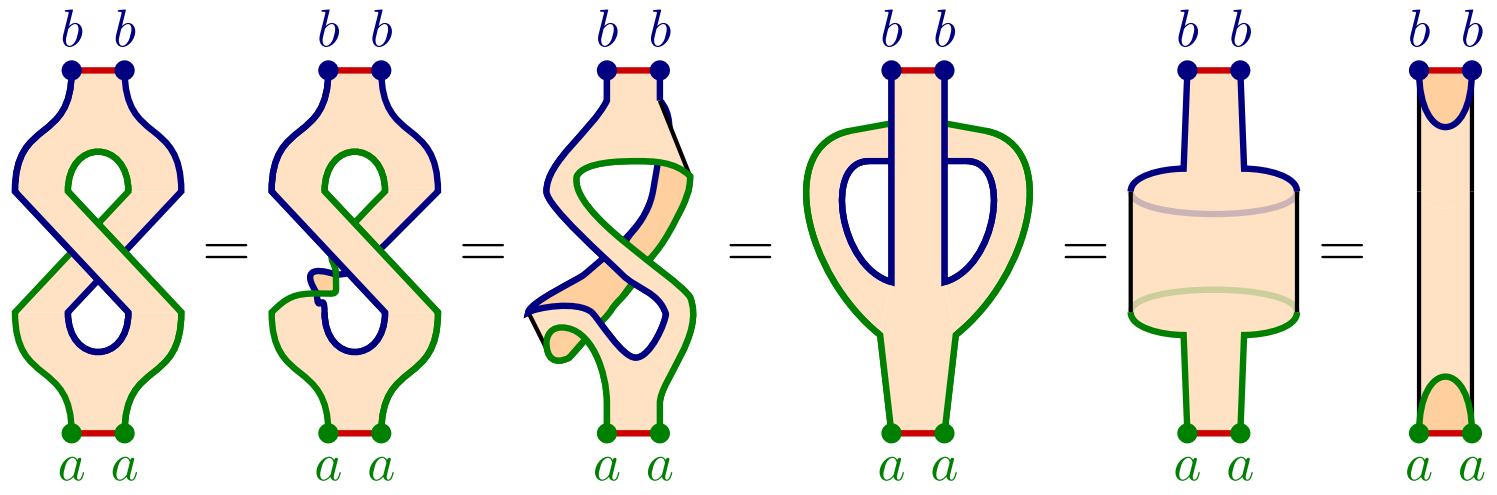

$$. \tag{3.13}$$

With the bordism category $\mathrm{Bord}_2^{\mathrm{oc}}(B)$ under explicit control, we now review how to equivalently encode the functor (3.3) in algebraic terms.

### 3.2.1 Closed sector

Restricting an open/closed TQFT $\mathcal{Z}^{\mathrm{oc}}$ as in (3.3) to the subcategory $\mathrm{Bord}_2 \subset \mathrm{Bord}_2^{\mathrm{oc}}(B)$ produces a closed TQFT. Hence $\mathcal{Z}^{\mathrm{oc}}(S^1)$ has the structure of a commutative Frobenius algebra $A^{\mathrm{c}}$. We write

$$\langle -,- \rangle^{\mathrm{c}} : A^{\mathrm{c}} \otimes_{\Bbbk} A^{\mathrm{c}} \longrightarrow \Bbbk \tag{3.14}$$

for its nondegenerate pairing.

### 3.2.2 Open sector

Let $\mathrm{Bord}_2^{\mathrm{o}}(B)$ be the subcategory of $\mathrm{Bord}_2^{\mathrm{oc}}(B)$ whose objects are only labelled intervals (and no circles), and whose morphisms are generated by the first five bordisms in (3.4), subject to (3.5)–(3.8) and the twist relations. By definition an *open TQFT* is a symmetric monoidal functor $\mathrm{Bord}_2^{\mathrm{o}}(B) \to \mathrm{Vect}_{\Bbbk}$.

We construct a category $\mathcal{C}^{\mathrm{o}}$ from $\mathcal{Z}^{\mathrm{oc}}$ restricted to $\mathrm{Bord}_2^{\mathrm{o}}(B)$ as follows. The set of objects is the set of boundary conditions,

$$\mathrm{Obj}(\mathcal{C}^{\mathrm{o}}) = B\,,$$

and Hom spaces are what $\mathcal{Z}^{\mathrm{oc}}$ assigns to labelled intervals,

$$\mathrm{Hom}(a,b) = \mathcal{Z}^{\mathrm{oc}}\Big($$ 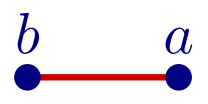 $$\Big).$$

Composition is defined to be

$$\mathcal{Z}^{\mathrm{oc}}\Bigg($$ 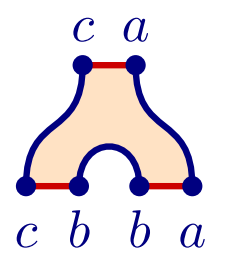 $$\Bigg) : \mathrm{Hom}(b,c) \times \mathrm{Hom}(a,b) \longrightarrow \mathrm{Hom}(a,c)\,.$$

It is associative by relation (3.5) and unital thanks to (3.6). This establishes that $\mathcal{C}^{\mathrm{o}}$ is a $\Bbbk$-linear category.

The relations (3.7) and (3.8) endow $\mathcal{C}^{\mathrm{o}}$ with additional structure. For all $a, b \in B$ there are $\Bbbk$-linear pairings

$$\langle -,-\rangle_{ab} = \mathcal{Z}^{\mathrm{oc}}\Big($$
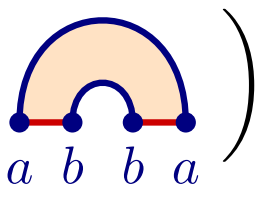

$$\Big) : \mathrm{Hom}(b,a) \otimes_{\Bbbk} \mathrm{Hom}(a,b) \longrightarrow \Bbbk \,. \qquad (3.15)$$

By relation (3.8) these pairings are symmetric in the sense that $\langle \Phi, \Psi\rangle_{ab} = \langle \Psi, \Phi\rangle_{ba}$, and they are nondegenerate since according to (3.6) and (3.7). Finally, the pairings are compatible with multiplication, $\langle \Phi_1, \Phi_2\Phi_3\rangle_{ab} = \langle \Phi_1\Phi_2, \Phi_3\rangle_{ac}$, by definition of the pairing and associativity (3.5).

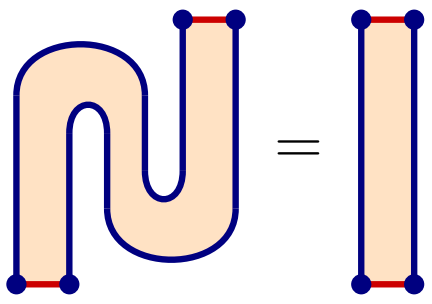

It follows from the above that the vector spaces $\mathrm{End}(a)$ have the structure of a symmetric Frobenius algebra. Hence one can think of $\mathcal{C}^{\mathrm{o}}$ as 'many Frobenius algebras glued together'. Unfortunately, the name 'Frobenius category' was already taken,[10] and instead the term 'Calabi-Yau' category is used.

To give the definition we first broaden the context. A *Serre functor* on a $\Bbbk$-linear category $\mathcal{C}$ is a functor $\Sigma : \mathcal{C} \to \mathcal{C}$ together with isomorphisms

$$\eta_{ab} : \mathrm{Hom}\big(a,b\big) \xrightarrow{\ \cong\ } \mathrm{Hom}\big(b, \Sigma(a)\big)^*$$

which are natural in $a, b \in \mathcal{C}$. From this one obtains the nondegenerate *Serre pairings*

$$\langle -,-\rangle_{ab} : \mathrm{Hom}\big(b, \Sigma(a)\big) \otimes_{\Bbbk} \mathrm{Hom}\big(a,b\big) \longrightarrow \Bbbk \,, \qquad \Psi \otimes \Phi \longmapsto \eta_{aa}(1_a)(\Psi\Phi)$$

by duality. By definition a *Calabi-Yau category* is a $\Bbbk$-linear category $\mathcal{C}$ together with a trivial Serre functor $\Sigma = 1_{\mathcal{C}}$. (If $\mathcal{C}$ is a *triangulated* category, then the Serre functor may be the identity only up to a shift.)

The eponymous example of a (triangulated) Calabi-Yau category is the bounded derived category of coherent sheaves on a Calabi-Yau variety. Indeed, as was noted in Section 2.4.2, the Serre functor on $\mathbb{D}^{\mathrm{b}}(\mathrm{coh}(U))$ for any smooth and proper variety $U$ is $\Sigma_U \cong \omega_U[\dim_{\mathbb{C}} U] \otimes_{\mathbb{C}} (-)$. But $U$ is Calabi-Yau iff the canonical line bundle $\omega_U$ is trivial. In light of the results of Section 3.3 below, Section 2.4 provides many further examples of Calabi-Yau categories.

It follows that the category $\mathcal{C}^{\mathrm{o}}$ we constructed from the open TQFT $\mathcal{Z}^{\mathrm{oc}}|_{\mathrm{Bord}_2^{\mathrm{o}}(B)}$ is Calabi-Yau with Serre pairings (3.15). It was shown in [Laz, MS, LP1] that the converse is also true:

[10] A *Frobenius category* is a Quillen exact category with enough injectives and enough projectives, such that injectives and projectives coincide.

**Theorem 3.3.** The above construction is an equivalence of groupoids between open TQFTs and Calabi-Yau categories.

### 3.2.3 Open/closed sector

It remains to work out the algebraic meaning of the last two generators in (3.4) as well as their relations (3.9)–(3.12). For this, we first define the *bulk-boundary maps* to be

$$\beta_a := \mathcal{Z}^{\mathrm{oc}}\Bigg($$
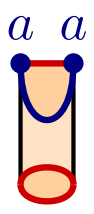
$$\Bigg) : A^{\mathrm{c}} \longrightarrow \mathrm{End}(a)$$

for all $a \in B$. It maps the closed sector (or 'bulk theory') to the open sector (with 'boundary condition' $a$).

Due to relations (3.9) and (3.10), $\beta_a$ is a map of algebras into the centre of $\mathrm{End}(a)$, i. e.

$$\beta_a(\phi)\,\Psi = \Psi\,\beta_a(\phi)$$

for all $\phi \in A^{\mathrm{c}}$ and $\Psi \in \mathrm{End}(a)$. Furthermore, the *boundary-bulk map* (which need not be a map of algebras)

$$\beta^a := \mathcal{Z}^{\mathrm{oc}}\Bigg($$
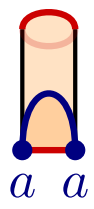
$$\Bigg) : \mathrm{End}(a) \longrightarrow A^{\mathrm{c}}$$

is adjoint to $\beta_a$ with respect to the Frobenius pairings by (3.11):

$$\Big\langle \beta^a(\Psi), \phi \Big\rangle^{\mathrm{c}} = \Big\langle \Psi, \beta_a(\phi) \Big\rangle_{aa} .$$

The most interesting condition on the maps $\beta^a$ derives from (3.12). For $\Phi \in \mathrm{End}(a)$ and $\Psi \in \mathrm{End}(b)$ let us consider the map

$${}_{\Psi}m_{\Phi} : \mathrm{Hom}(a,b) \longrightarrow \mathrm{Hom}(a,b)\,, \quad \Omega \longmapsto \Psi\Omega\Phi\,.$$

Note that for $\Phi = \Psi = 1_a \in \mathrm{End}(a)$, the map ${}_{1_a}m_{1_a}$ is simply the identity operator on $\mathrm{End}(a)$. Then for any open/closed TQFT $\mathcal{Z}^{\mathrm{oc}}$ as above we have:

**Theorem 3.4** (Cardy condition)**.** Let $\Phi \in \mathrm{End}(a)$ and $\Psi \in \mathrm{End}(b)$. Then

$$\mathrm{tr}\left({}_{\Psi}m_{\Phi}\right) = \Big\langle \beta^b(\Psi), \beta^a(\Phi) \Big\rangle^{\mathrm{c}} . \tag{3.16}$$

The significance of the Cardy condition is that a trace in the open sector can be computed from a pairing in the closed sector. The qualifier 'theorem' is more than appropriate: as shown in [CW], in the special case of the Calabi-Yau category being of the form $\mathbb{D}^{\mathrm{b}}(\mathrm{coh}(U))$, the Cardy condition is the Hirzebruch-Riemann-Roch theorem!

To prove Theorem 3.4 we use duality to rewrite (3.13) as

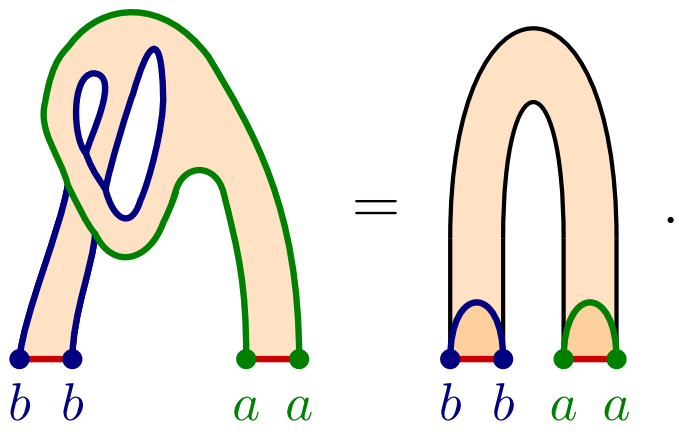

(3.17)

Let us choose bases $\{e_i\}$ and $\{\widetilde{e}_i\}$ of $\mathrm{Hom}(a,b)$ and $\mathrm{Hom}(b,a)$, respectively, such that we have $\mathcal{Z}^{\mathrm{oc}}($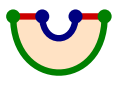$)(1) = \sum_i \widetilde{e}_i \otimes e_i$ for the copairing. Then thanks to (3.6) and (3.7) the basis $\{e_i^* = \langle \widetilde{e}_i, - \rangle_{ab}\}$ is dual to $\{e_i\}$, and we compute

$$\begin{aligned}
\mathrm{tr}\left({}_\Psi m_\Phi\right) &= \sum_i e_i^*\left(\Psi e_i \Phi\right) \\
&= \sum_i \left\langle \widetilde{e}_i, \Psi e_i \Phi \right\rangle_{ab} \\
&= \mathcal{Z}^{\mathrm{oc}}\Big(\text{LHS of } (3.17)\Big)(\Psi \otimes \Phi) \\
&= \mathcal{Z}^{\mathrm{oc}}\Big(\text{RHS of } (3.17)\Big)(\Psi \otimes \Phi) \\
&= \left\langle \beta^b(\Psi), \beta^a(\Phi) \right\rangle^{\mathrm{c}}.
\end{aligned}$$

In summary, a 2-dimensional open/closed TQFT has the following algebraic description, generalising the fact that closed TQFTs are equivalent to commutative Frobenius algebras:

**Theorem 3.5.** The construction reviewed in this section gives an equivalence between open/closed TQFTs and the following data:

- a commutative Frobenius algebra $A^{\mathrm{c}}$ with nondegenerate pairing $\langle -,- \rangle^{\mathrm{c}}$,
- a Calabi-Yau category $\mathcal{C}^{\mathrm{o}}$,
- $\Bbbk$-linear maps $\beta_a : A^{\mathrm{c}} \rightleftarrows \mathrm{End}(a) : \beta^a$ for all $a \in \mathcal{C}^{\mathrm{o}}$,

such that

(i) $\beta_a$ are algebra maps with image in the centre,

(ii) $\beta_a$ and $\beta^a$ are adjoint with respect to the Frobenius pairings,

(iii) the Cardy condition

$$\mathrm{tr}\left({}_\Psi m_\Phi\right) = \left\langle \beta^b(\Psi), \beta^a(\Phi) \right\rangle^{\mathrm{c}}$$

holds for all $\Phi \in \mathrm{End}(a)$ and $\Psi \in \mathrm{End}(b)$ in $\mathcal{C}^{\mathrm{o}}$.

## 3.3 Open/closed TQFTs from pivotal 2-categories

In Section 3.1 we constructed a closed TQFT for every object $\alpha$ in a bicategory $\mathcal{B}$ satisfying Assumption 3.1. Now we complete the construction by assigning an open/closed TQFT to every $\alpha \in \mathcal{B}$.

### 3.3.1 Open sector

The natural starting point for the open sector is the category

$$\mathcal{C}_\alpha = \mathcal{B}(\mathbb{O}, \alpha) \tag{3.18}$$

of 1-morphisms between $\alpha$ and the 'trivial theory', i. e. the monoidal unit $\mathbb{O} \in \mathcal{B}$. In the interpretation of 1-morphisms as defect lines it is immediate to view them as boundary conditions:

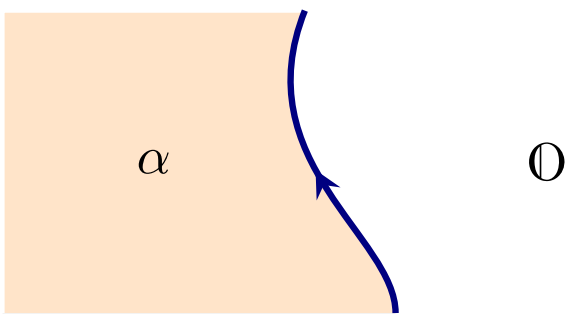

We claim that (3.18) comes with the structure of a Calabi-Yau category. It certainly is $\Bbbk$-linear, but where are the nondegenerate pairings? They are encoded in the duality structure of $\mathcal{B}$ together with the assumption $\mathcal{B}(\mathbb{O}, \mathbb{O}) \cong \mathrm{Vect}_\Bbbk$. Indeed, for $X, Y \in \mathcal{B}(\mathbb{O}, \alpha)$ we have

$$X^\dagger \otimes Y \in \mathrm{Vect}_\Bbbk \,,$$

hence the Zorro moves for $X^\dagger \otimes Y$ tell us that

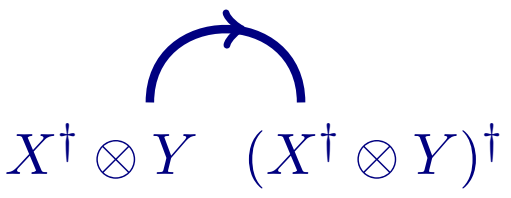

$$: \big(X^\dagger \otimes Y\big) \otimes_\Bbbk \big(X^\dagger \otimes Y\big)^\dagger \longrightarrow \Bbbk$$

is a nondegenerate pairing. Precomposing with the isomorphism

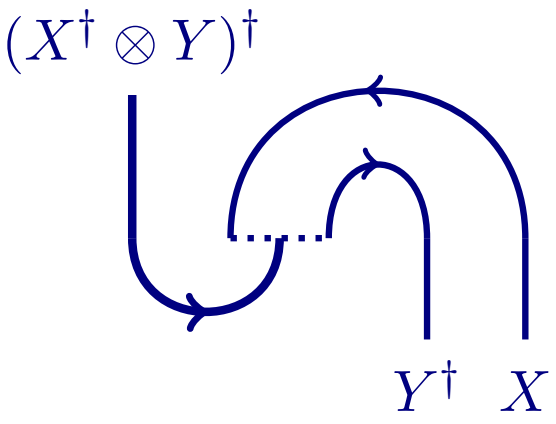

$$: Y^\dagger \otimes X \longrightarrow \big(X^\dagger \otimes Y\big)^\dagger$$

and using the Zorro move for $X^\dagger \otimes Y$ again, we obtain the nondegenerate pairing

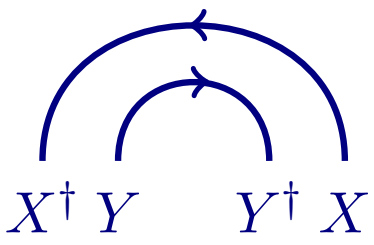

$$: \big(X^\dagger \otimes Y\big) \otimes_\Bbbk \big(Y^\dagger \otimes X\big) \longrightarrow \Bbbk \,. \tag{3.19}$$

To translate the above into a pairing in $\mathcal{C}_\alpha$, we use the isomorphism of vector spaces

$$\mathrm{Hom}(X,Y) \xrightarrow{\cong} X^\dagger \otimes Y\,, \qquad (3.20)$$

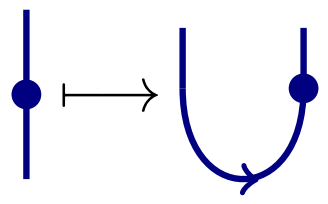

with inverse

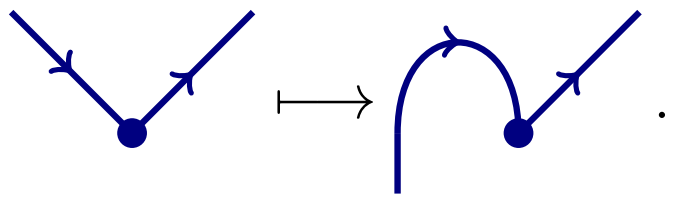

Precomposing (3.19) with (3.20) we obtain the Serre pairing

$$\langle -,-\rangle^\alpha_{YX} : \mathrm{Hom}(X,Y) \otimes_{\Bbbk} \mathrm{Hom}(Y,X) \longrightarrow \Bbbk$$

with

$$\langle \Phi, \Psi\rangle^\alpha_{YX} =$$

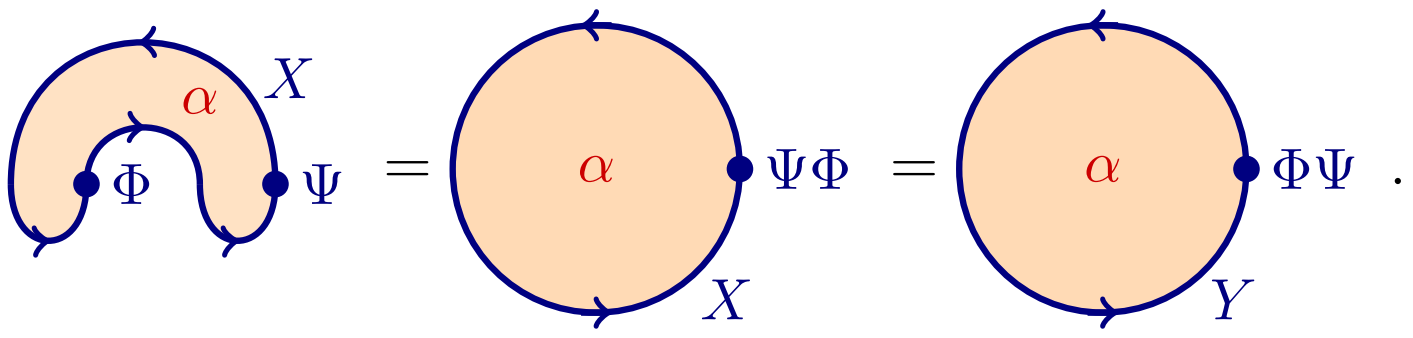

Since this is the left trace (recall (2.9)) of $\Psi\Phi \in \mathrm{End}(X)$, the pairings $\langle -,-\rangle^\alpha_{YX}$ are symmetric: $\langle \Phi,\Psi\rangle^\alpha_{YX} = \langle \Psi,\Phi\rangle^\alpha_{XY}$.

In summary, we have proved

**Theorem 3.6.** Let $\mathcal{B}$ be a bicategory satisfying Assumption 3.1. Then for every $\alpha \in \mathcal{B}$, the category $\mathcal{B}(\mathbb{O}, \alpha)$ naturally has the structure of a Calabi-Yau category with pairing

$$\langle -,-\rangle^\alpha_{YX} : \mathrm{Hom}(X,Y) \otimes_{\Bbbk} \mathrm{Hom}(Y,X) \longrightarrow \Bbbk\,, \quad \langle \Phi, \Psi\rangle^\alpha_{YX} =$$

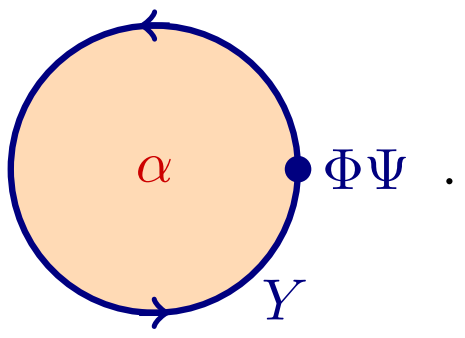

Recall that for the closed sector $A_\alpha = \mathrm{End}(1_\alpha)$ in Section 3.1, the nondegenerate pairing (3.2) was the 'sphere correlator projected to a disc'. But in the open sector the role of the sphere correlator is played by the disc correlator – so we could have guessed the result of Theorem 3.6 from the start!

### 3.3.2 Open/closed sector

We continue to work with a bicategory $\mathcal{B}$ satisfying Assumption 3.1. For every $\alpha \in \mathcal{B}$ we have already obtained a commutative Frobenius algebra $A_\alpha = \mathrm{End}(1_\alpha)$ and a Calabi-Yau category $\mathcal{C}_\alpha = \mathcal{B}(\mathbb{O}, \alpha)$. According to Theorem 3.5, we also need maps

$$\beta_X : A_\alpha \rightleftarrows \mathrm{End}(X) : \beta^X$$

for all $X \in \mathcal{C}_\alpha$ in order to associate an open/closed TQFT to $\alpha$. In the graphical calculus, intuition about these maps becomes a rigorous definition: we set the bulk-boundary map to be

$$\beta_X : \mathrm{End}(1_\alpha) \longrightarrow \mathrm{End}(X)\,, \tag{3.21}$$

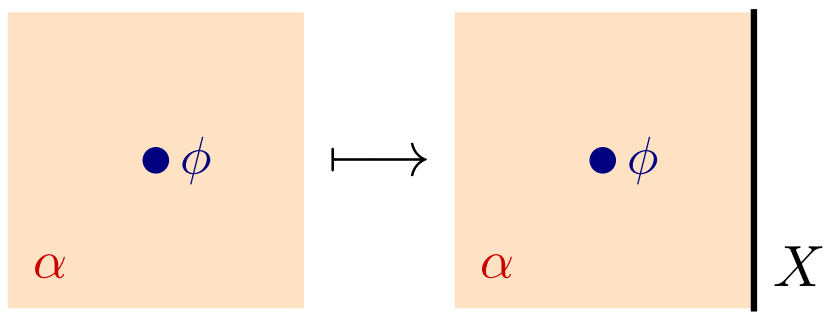

and the boundary-bulk map is

$$\beta^X : \mathrm{End}(X) \longrightarrow \mathrm{End}(1_\alpha)\,, \tag{3.22}$$

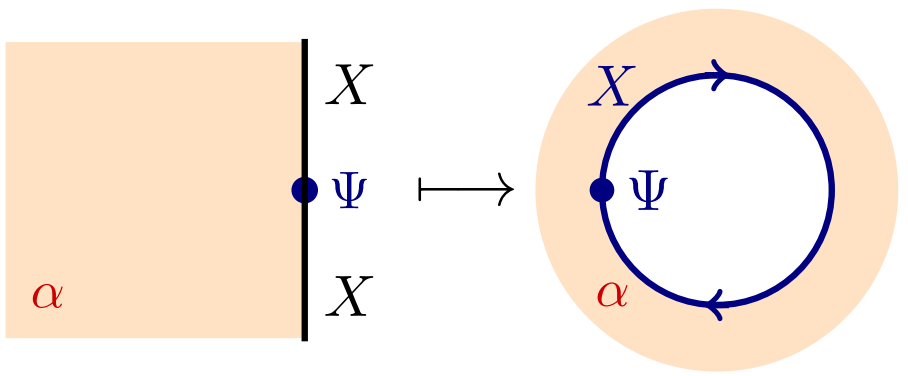

We have to check that the maps $\beta_X, \beta^X$ satisfy the three conditions (i)–(iii) in Theorem 3.5. The first one is easy: $\beta_X$ is manifestly an algebra map into the centre, thanks to the coherence theorem behind the graphical calculus.

To prove conditions (ii) and (iii), we make one additional assumption. For motivation, recall that we obtained the nondegenerate pairing $\langle -,- \rangle_\alpha$ on $A_\alpha$ in (3.2) by projecting an $\alpha$-decorated sphere onto a plane. Now we consider a sphere with two incoming boundary circles cut out, together with a $D$-decorated defect line which separates two phases labelled $\alpha$ and $\beta$ as follows:

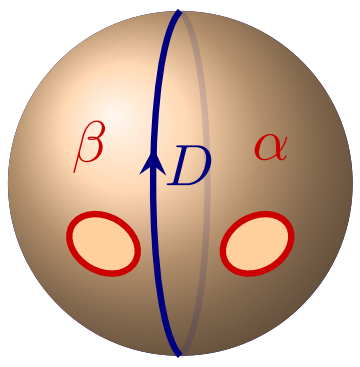

By isotopically deforming the defect line we obtain the identity

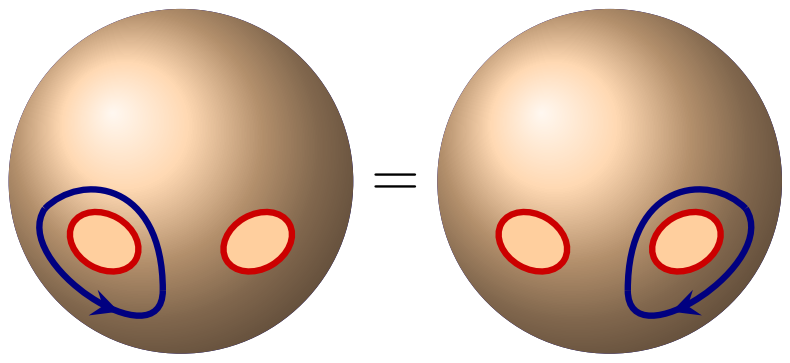

in the defect bordism category. Hence we expect the following property for the bicategory $\mathcal{B} = \mathcal{B}_\mathcal{Z}$ associated to a defect TQFT:

**Assumption 3.7.** For every $D \in \mathcal{B}(\alpha,\beta)$ and $\Psi \in \mathrm{End}(D)$, the pairings $\langle -,- \rangle_\alpha$ on $\mathrm{End}(1_\alpha)$ satisfy

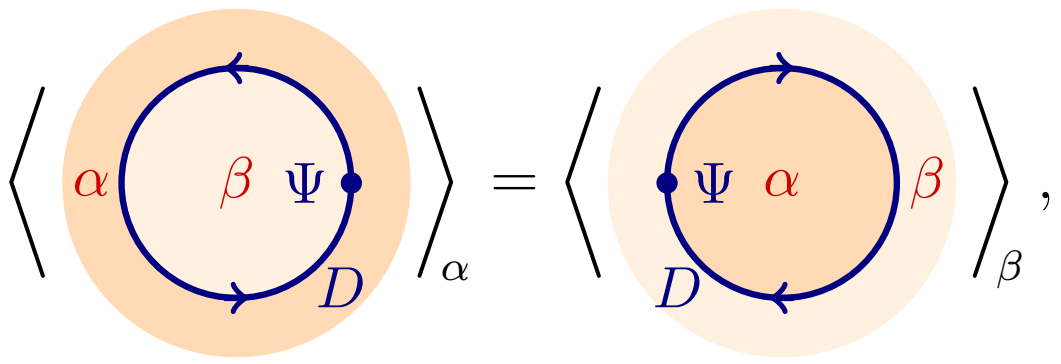

(3.23)

where we write $\langle \phi \rangle_\alpha$ for $\langle \phi, 1 \rangle_\alpha$.

Under this assumption we can verify condition (ii) of Theorem 3.5 in just one line. Let $\phi \in \mathrm{End}(1_\alpha)$, $X \in \mathcal{B}(\mathbb{O}, \alpha)$ and $\Psi \in \mathrm{End}(X)$. Then

$$\Big\langle \beta_X(\phi), \Psi \Big\rangle_{XX} =$$
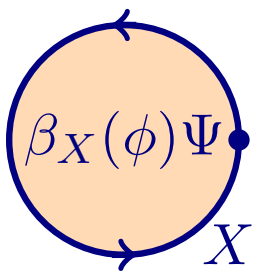

$$= \Big\langle$$
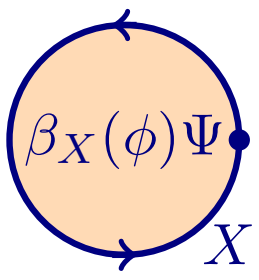

$$\Big\rangle_{\mathbb{O}} = \Big\langle$$
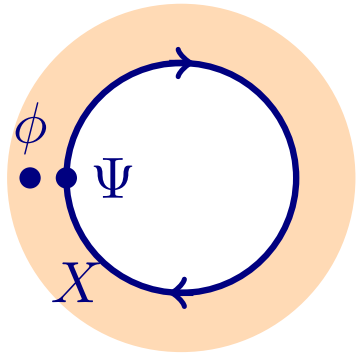

$$\Big\rangle_\alpha = \Big\langle \phi, \beta^X(\Psi) \Big\rangle^{\mathrm{c}}$$

so indeed $\beta_X$ and $\beta^X$ are adjoint with respect to one another.

At last we prove the Cardy condition. Let $X, Y \in \mathcal{B}(\mathbb{O}, \alpha)$, $\Phi \in \mathrm{End}(X)$ and $\Psi \in \mathrm{End}(Y)$. Then

$$\Big\langle \beta^X(\Phi), \beta^Y(\Psi) \Big\rangle_\alpha \overset{(3.22)}{=} \Big\langle$$
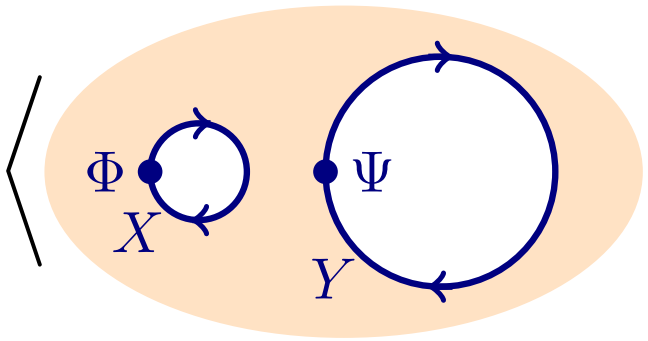

$$\Big\rangle_\alpha$$

$$\overset{(3.23)}{=} \Big\langle$$
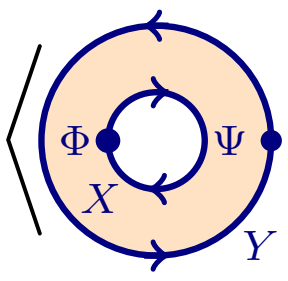

$$\Big\rangle_{\mathbb{O}} =$$
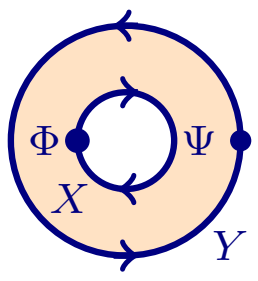

$$\overset{\text{Zorro}}{=}$$
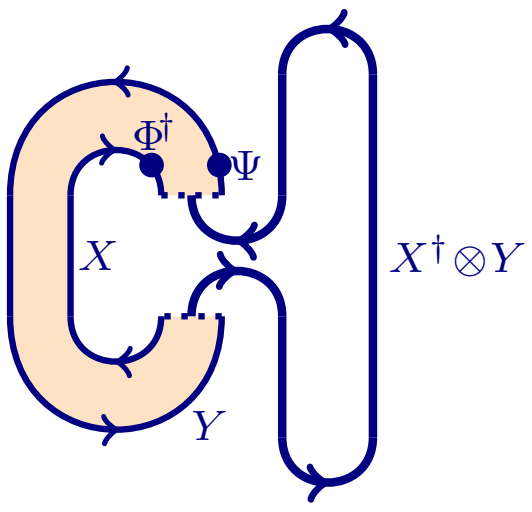

$$\overset{(2.7)}{=}$$
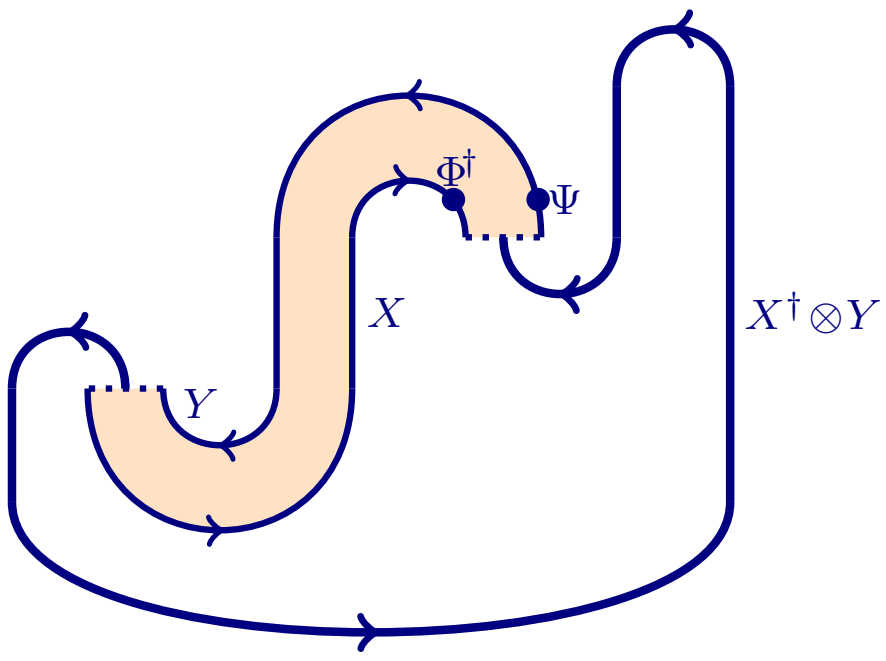

$$\overset{\text{Zorro}}{=}$$
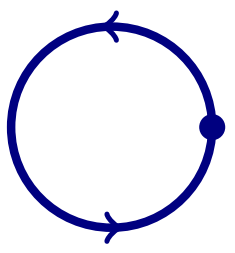
$$\Phi^\dagger \otimes \Psi \;= \mathrm{tr}\left( {}_\Psi m_\Phi \right).$$

To understand the last step, pick a basis $\{e_i\}$ of $X^\dagger \otimes Y \in \mathrm{Vect}_{\Bbbk}$, and write $\{e_i^*\}$ for the dual basis. The canonical isomorphisms $X^\dagger \otimes Y \cong \mathrm{Hom}(X, Y)$ and

$(X^\dagger \otimes Y)^* \cong \mathrm{Hom}(X,Y)^*$ act as

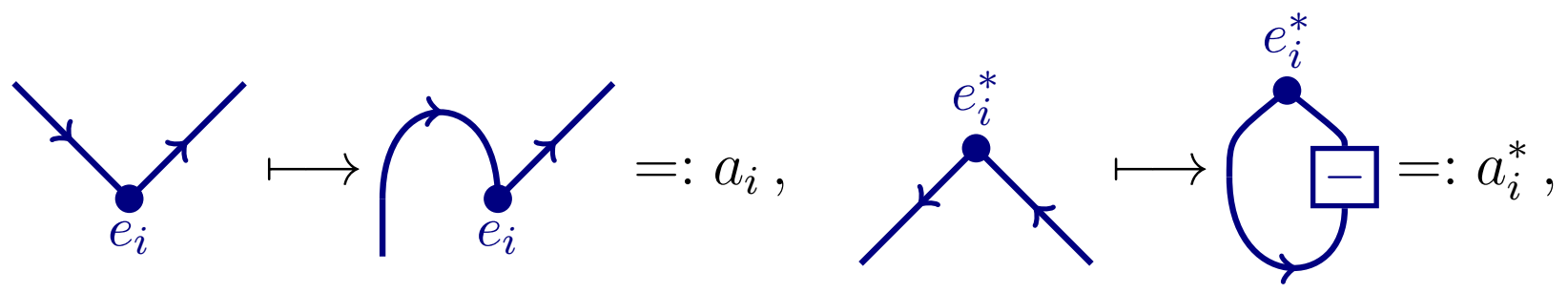

respectively. The Zorro move tells us that $a_j^*(a_i) = \delta_{ij}$, so $\{a_i^*\}$ is dual to the basis $\{a_i\}$ of $\mathrm{Hom}(X,Y)$. Hence we compute

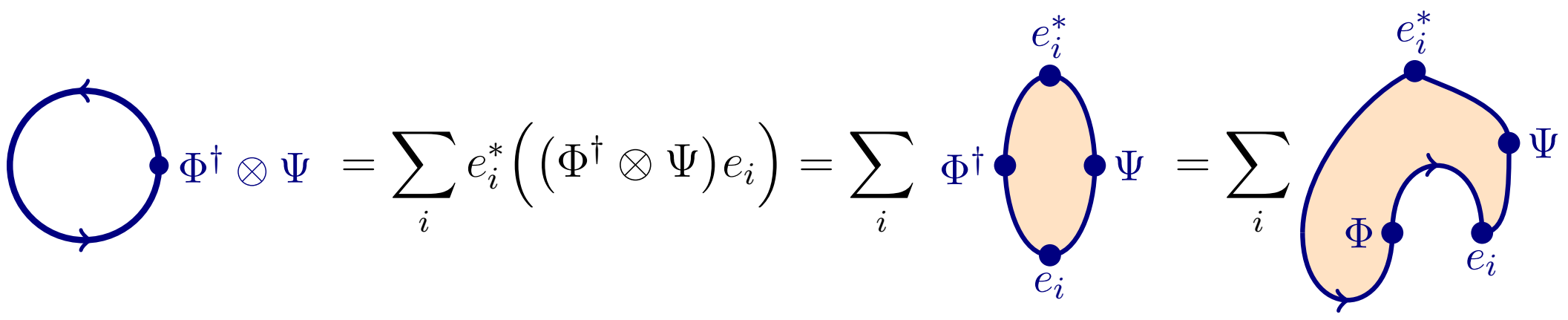

$$= \sum_i a_i^*\Big(\Psi a_i \Phi\Big) = \mathrm{tr}\left({}_\Psi m_\Phi\right).$$

This completes our construction of an open/closed TQFT from $\alpha$:

**Theorem 3.8.** Let $\mathcal{B}$ be a bicategory satisfying Assumptions 3.1 and 3.7. Then by the above construction for every $\alpha \in \mathcal{B}$ we have that

- $\mathrm{End}(1_\alpha)$ has the structure of a commutative Frobenius algebra,
- $\mathcal{B}(\mathbb{O},\alpha)$ has the structure of a Calabi-Yau category,
- for every $X \in \mathcal{B}(\mathbb{O},\alpha)$, the maps $\beta_X : \mathrm{End}(1_\alpha) \rightleftarrows \mathrm{End}(X) : \beta^X$ of (3.21) and (3.22) satisfy the conditions in Theorem 3.5.